\documentclass[11pt]{article}
\usepackage[titletoc]{appendix}
\usepackage[labelsep=period]{caption}
\usepackage{bm}
\usepackage{fullpage}
\usepackage{epsfig}
\usepackage{graphicx}
\usepackage{subcaption}
\usepackage{latexsym}
\usepackage{amsmath}
\usepackage{amsfonts}
\usepackage{amssymb}
\usepackage{mathrsfs}
\usepackage{pifont}
\usepackage{yhmath}
\usepackage{bbm}
\usepackage{dsfont}
\usepackage{eufrak}
\usepackage{wasysym}
\usepackage{underscore}
\usepackage{epstopdf}
\usepackage{fontenc}
\usepackage{amsthm}
\newtheorem{theorem}{Theorem}[section]

\newtheorem{corollary}[theorem]{Corollary}

\newtheorem{lemma}[theorem]{Lemma}

\def\qed{\hfill \rule{4pt}{7pt}}
\usepackage{fullpage,graphicx}
\newcommand{\de}{\backslash}
\DeclareMathAlphabet{\mathpzc}{OT1}{pzc}{m}{it}
\title{\bf Integral Biflow Maximization}
\author{\vspace{2mm}  Guoli Ding$^{a}$
\quad Rongchuan Tao$^{b}$
\quad Mengxi Yang$^{b}$\thanks{Corresponding author. E-mail: yangmx@connect.hku.hk.}
\quad Wenan Zang$^{b}$\thanks{Supported in part by the Research Grants Council of Hong Kong.} \\
$\stackrel{a}{}$ Mathematics Department, Louisiana State  University\\ Baton Rouge, LA 70803, USA \smallskip\\
$\stackrel{b}{}$ Department of Mathematics, University of Hong Kong \\ Hong Kong, China}
\begin{document}
\date{}
\maketitle

\begin{abstract}

Let $G=(V,E)$ be a graph with four distinguished vertices, two sources $s_1, s_2$ and two sinks $t_1,t_2$, let $c:\, E \rightarrow \mathbb Z_+$ be a 
capacity function, and let ${\cal P}$ be the set of all simple paths in $G$ from $s_1$ to $t_1$ or from $s_2$ to $t_2$. A {\em biflow} (or {\em $2$-commodity flow}) 
in $G$ is an assignment $f:\, {\cal P}\rightarrow \mathbb R_+$ such that $\sum_{e \in Q \in {\cal P}}\, f(Q) \le c(e)$ for all $e \in E$, 
whose {\em value} is defined to be $\sum_{Q \in {\cal P}}\, f(Q)$. A {\em bicut} in $G$ is a subset $K$ of $E$ that contains at least one edge from each 
member of ${\cal P}$, whose {\em capacity} is $\sum_{e\in K}\, c(e)$. In 1977 Seymour characterized, in terms of forbidden structures, 
all graphs $G$ for which the max-biflow (integral) min-bicut theorem holds true (that is, the maximum value of an integral biflow is equal to 
the minimum capacity of a bicut for every capacity function $c$); such a graph $G$ is referred to as a Seymour graph. Nevertheless, his proof is not 
algorithmic in nature. In this paper we present a combinatorial polynomial-time algorithm for finding maximum integral biflows in Seymour graphs, which 
relies heavily on a structural description of such graphs.  

\vskip 4mm

\noindent {\bf MSC 2000 subject classification.} Primary:  90C27, 05C85, 68Q25.

\noindent {\bf OR/MS subject classification.} Primary: Programming/graphs.

\noindent {\bf Key words.} Biflow, bicut, algorithm, structure, characterization.

\end{abstract}

\newpage

\section{Introduction}

Due to its deep intellectual content and wide range of applicability, network flow theory has attracted tremendous 
research efforts over the past seven decades, and has become a flourishing discipline containing a body of beautiful 
theorems and powerful methods.  One of the most fundamental problems in network flow theory is the {\em maximum flow problem}, 
which aims to ship as much of a single commodity as possible from the source to the sink in a given network, without 
exceeding the capacity of each arc. As is well known, this problem can be solved in strongly polynomial time.  Furthermore, 
its solution can be guaranteed to be integral provided all capacities are integers. Among various approaches to the 
maximum flow problem, the augmenting path algorithm obviously possesses the greatest mathematical beauty, because 
it not only finds a solution to the dual, the {\em minimum cut problem}, as a by-product, but also yields a proof 
of the celebrated max-flow min-cut theorem by Ford and Fulkerson \cite{FF}, asserting that the maximum value of a flow from the source 
to the sink in a capacitated network $N$ equals the minimum capacity of a cut, regardless of the structure of $N$. This central theorem 
in network flow theory permits us to model a rich variety of problems as maximum flow or minimum cut problems, even though 
on the surface lots of them do not appear to have a network flow structure. So it has far-reaching implications.

In many other problem settings, two commodities, each governed by its own flow constraints, share the same network and common 
capacities. To find a desired biflow (or 2-commodity flow), we have to solve the individual single commodity problems in concert 
with each other, due to the dependence between them. Now a natural question to ask is: Can we generalize the aforementioned 
integrality, duality, and polynomial solvability properties satisfied by the maximum flow problem to the 2-commodity case?  This
question on undirected networks has been a subject of extensive research; it is also the major concern of the present paper.

We introduce some notation before proceeding. As usual, $\mathbb Z_+$ stands for the set of nonnegative integers and $\mathbb R_+$ 
for the set of nonnegative real numbers. Let $I$ be a set and let $\alpha$ be a real-valued function with domain $I$. For each 
subset $S$ of $I$, let $\alpha(S):=\sum_{s\in S} \,\alpha(s)$. For each vertex $v$ of a graph $G$, we use $\delta_G(v)$ 
to denote the set of all edges incident with $v$;  we shall drop the subscript $G$ if there is no danger of
confusion.  

Let $G=(V,E)$ be a graph with four distinguished vertices, two sources $s_1, s_2$ and two sinks $t_1,t_2$, 
let $c:\, E \rightarrow \mathbb Z_+$ be a capacity function, and let ${\cal P}$ be the set of all simple paths in $G$ from $s_1$ to $t_1$ 
or from $s_2$ to $t_2$. A {\em biflow} (or {\em $2$-commodity flow}) in $G$ is an assignment $f:\, {\cal P}\rightarrow \mathbb R_+$ 
such that $\sum_{e \in Q \in {\cal P}}\, f(Q) \le c(e)$ for all $e \in E$, whose {\em value} is defined to be $\sum_{Q \in {\cal P}}\, f(Q)$. 
The {\em maximum biflow problem} is to find a biflow in $G$ with maximum value, and the {\em maximum integral biflow problem} is to 
find an integral biflow in $G$ with maximum value. Just as in the one commodity case, cut also plays an important role in the study of 
these problems. A subset $K$ of $E$ is called a {\em bicut} if it contains at least one edge from each member of ${\cal P}$; 
its {\em capacity} is defined to be $c(K)$. Hu \cite{Hu} proved that the maximum value of a biflow is equal to the minimum capacity of a 
bicut for every capacity function $c$. Furthermore, the maximum is achieved by a half-integral biflow. Sakarovitch \cite{S} and Itai \cite{I} devised different 
polynomial-time algorithms for finding such a biflow. Rothschild and Whinston \cite{RW} strengthened Hu's theorem \cite{Hu} by showing that if $c(\delta(v))$ 
is even for every vertex $v$, then the maximum is attained by an integral biflow, which can also be found in polynomial time. Rajagopalan \cite{R}
further extended Rothschild and Whinston's result to the case when $c(\delta(v))$ is even for every vertex $v \ne s_1,t_1,s_2,t_2$ (he
also showed a hole in the proof by Sakarovitch \cite{S} of this).  

In sharp contrast, the maximum integral biflow problem is $NP$-hard, as proved by Even, Itai and Shamir \cite{EIS}. For the case when
$G+\{s_1t_1,s_2t_2\}$ is planar, Lomonosov \cite{L} characterized for a {\em fixed} capacity function $c$, when the maximum value of an integral biflow and
the minimum capacity of a bicut are equal. He also showed that the maximum and the minimum differ by at most one for any capacity
function $c$.  Moreover, Korach and Penn \cite{KP} came up with a polynomial-time algorithm for the maximum integral biflow problem.  

In 1977 Seymour \cite{sey77} successfully characterized all matroids with the max-flow min-cut property; his theorem contains, as a corollary, a characterization 
of all graphs for which the max-biflow (integral) min-bicut theorem holds true; see also Theorem 71.2 in Schrijver \cite{s03}.  

\begin{figure}[h]
    \centering
    \includegraphics{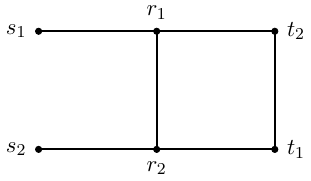}
    \caption{The forbidden subgraph $F$}
    \label{fig:forbidden structure}
\end{figure}

\begin{theorem}[Max-biflow (integral) Min-bicut Theorem; Seymour \cite{sey77}] \label{Thm1} 
Let $G=(V,E)$ be a graph with two sources $s_1, s_2$ and two sinks $t_1,t_2$. Then the maximum value of an integral biflow is equal to the minimum capacity 
of a bicut for every capacity function $c:\, E \rightarrow \mathbb Z_+$, if and only if $G$ contains no subgraph contractible to the graph $F$ depicted in 
Figure \ref{fig:forbidden structure}, up to interchanging $s_1$ and $t_1$, and $s_2$ and $t_2$. 
\end{theorem}

So contractions of $F$ are the only obstructions to the desired min-max relation, Seymour's proof, however, is not algorithmic in nature. For simplicity, we refer
to any good graph exhibited in the above theorem as a {\em Seymour graph}. The purpose of this paper is to present a combinatorial polynomial-time 
algorithm for finding maximum integral biflows in Seymour graphs, which relies heavily on a structural description of such graphs to be given in
the next section.  

\begin{theorem}\label{Thm2} 
Let $G=(V,E)$ be a Seymour graph with two sources $s_1, s_2$ and two sinks $t_1,t_2$ and let $c:\, E \rightarrow \mathbb Z_+$ be a capacity function.
Then a maximum integral biflow in $G$ can be found in $O(n^3+m^2)$ time, where $n=|V|$ and $m=|E|$.
\end{theorem}

The remainder of this paper is organized as follows. In Section 2, we present a structural description of Seymour graphs, which asserts that
Seymour graphs essentially fall into four classes. In Section 3, we give a brief review of network flow theory, and prove a technical lemma, 
which will be used repeatedly in the subsequent sections. In Sections 4-6, we devise a combinatorial polynomial-time algorithm for finding
maximum integral biflows in these four classes of Seymour graphs, respectively.

\section{Structural Description}

We make some preparations before presenting a structural description of Seymour graphs and a proof of its validity.

Let $G$ be a graph. We use $V(G)$ and $E(G)$ to denote its vertex set and edge set, respectively. For each $A \subseteq V(G)$, we use 
$G\de A$ to denote the graph obtained from $G$ by deleting all vertices in $A$, and write $G\de a$ for $G\de A$ when $A=\{a\}$. Similarly,
for each $B \subseteq E(G)$, we use $G\de B$ to denote the graph obtained from $G$ by deleting all edges in $B$, and write $G\de b$ 
for $G\de B$ when $B=\{b\}$. 

Given two graphs $G_1=(V_1,E_1)$ and $G_2=(V_2,E_2)$, let $G_1\cup G_2:=(V_1\cup V_2, E_1\cup E_2)$ and $G_1\cap G_2:=(V_1\cap V_2, E_1 \cap E_2)$; we 
call $G_1\cup G_2$ and $G_1\cap G_2$ the {\em union} and {\em intersection} of $G_1$ and $G_2$, respectively. When $V_1\cap V_2=\emptyset$, a $3$-{\em sum} 
of $G_1$ and $G_2$ is obtained by first choosing a triangle (a cycle of length $3$) in each of $G_1$ and $G_2$, then identifying these two triangles, and 
finally deleting $0,1,2$, or $3$ of these identified edges. 

A {\it separation} of a graph $G$ is a pair $(G_1,G_2)$ of {\em edge-disjoint} non-spanning subgraphs of $G$ with $G_1 \cup G_2=G$; we call $V(G_1\cap G_2)$ 
the {\em separating set}, and call this set a $k$-{\em separating set} and this separation a $k$-{\em separation} if $|V(G_1\cap G_2)|=k$. We say that two 
edges $e_1$ and $e_2$ of $G$ are {\em separated} by the separation $(G_1,G_2)$ if each $G_i$ contains precisely one of $e_1$ and $e_2$ for $i=1,2$. We also 
say that two disjoint vertex subsets $U_1$ and $U_2$ of $G$ are {\em separated} by $(G_1,G_2)$ if each $G_i$ contains precisely one of $U_1$ and $U_2$ for $i=1,2$.

Let $H$ be a proper subgraph of a graph $G$. An {\it $H$-bridge} of $G$ is either a subgraph of $G$ induced by the edges of a component $C$ of $G\de V(H)$ together with 
the edges linking $C$ to $H$, or a subgraph induced by an edge not in $H$ but with both ends in $H$. We call the second type of bridges {\it trivial}, and
call the vertices of an $H$-bridge that are in $H$ its {\it feet}.

As usual, a graph $H$ is called a {\it minor} of a graph $G$ if $H$ arises from a subgraph of $G$ by contracting edges. We say that $G$ has a $H$-{\em minor}
if $G$ contains $H$ as a minor and that $G$ is {\it $H$-free} otherwise. 

\begin{figure}[h]
    \centering
    \begin{subfigure}{0.5\textwidth}
    \centering
    \includegraphics{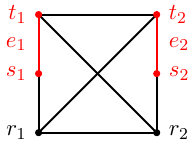} 
    \end{subfigure}
    \hspace{-3cm}
    \begin{subfigure}{0.5\textwidth}
    \centering
    \includegraphics{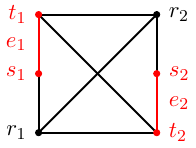}
    \end{subfigure}
    
    \caption{Two different drawings of $K_4^*$}
    \label{fig:K4*}
\end{figure}

Let $K_4^*$ denote the graph obtained from $F$ depicted in Figure \ref{fig:forbidden structure} by adding two edges $e_1=s_1t_1$ and $e_2=s_2t_2$  (up to interchanging 
$s_1$ and $t_1$, and $s_2$ and $t_2$); see Figure \ref{fig:K4*}. More generally, let $G$ be a graph with two specified pairs of vertices $\{s_1,t_1\}$ 
and $\{s_2,t_2\}$. The {\em augmented} graph of $G$, denoted by $\hat{G}$, is obtained from $G$ by adding an edge $e_i=s_it_i$ if $s_i$ and $t_i$ are nonadjacent 
for $i=1,2$. (So $K_4^*=\hat{F}$.) From Theorem \ref{Thm1}  we see that $G$ is a Seymour graph if and only if $\hat{G}$ is $K_4^*$-free. We shall demonstrate 
in Section 6 that the maximum integral biflow problem on a Seymour graph $G$ can be easily reduced to the case when $\hat{G}$ is $2$-connected and 
every $2$-separation of $\hat{G}$, if any, separates $e_1$ and $e_2$.

\begin{figure}[h]
    \centering
    \includegraphics{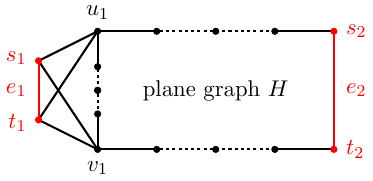}
    \caption{Graphs in $\mathcal{G}_2$}
    \label{fig:graph in G2}
\end{figure}

\begin{figure}[ht]
    \centering
    \begin{subfigure}{0.4\textwidth}
    \centering
    \includegraphics{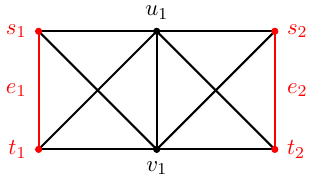} 
    \caption*{The graph in $\mathcal{G}_3$ with $H=K_2$}
    \end{subfigure}
    \begin{subfigure}{0.5\textwidth}
    \centering
    \includegraphics{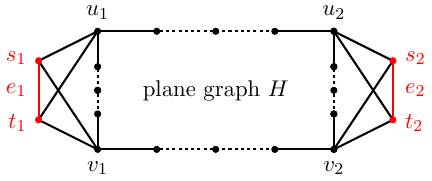}
    \caption*{Graphs in $\mathcal{G}_3$ with $H\neq K_2$}
    \end{subfigure}
    
    \caption{Graphs in $\mathcal{G}_3$}
    \label{fig:graph in G3}
\end{figure}

The following three classes of graphs with two specified nonadjacent edges $e_1=s_1t_1$ and $e_2=s_2t_2$ will be involved in our structural description:
\begin{itemize}
\vspace{-2mm}
\item $\mathcal{G}_1$ consists of all $2$-connected plane graphs $H$ with $e_1$ and $e_2$ on its outer facial cycle; 
\vspace{-2mm}
\item $\mathcal{G}_2$ consists of all graphs obtained from a $2$-connected plane graph $H$ with four distinct vertices $u_1,v_1, s_2,t_2$ (together with edge $e_2$) on its 
outer facial cycle by adding two vertices $s_1, t_1$ and five edges $s_1t_1, s_1u_1, s_1v_1, t_1u_1, t_1v_1$, up to interchanging $e_1$ and $e_2$ (thereby
interchanging $s_1$ and $s_2$ and interchanging $t_1$ and $t_2$ as well; see Figure \ref{fig:graph in G2}); and 
\vspace{-2mm}
\item $\mathcal{G}_3$ consists of all graphs obtained from the complete graph $H$ on two vertices $u_1, v_1$ (now define $u_2=u_1$ and $v_2=v_1$) or from a $2$-connected 
plane graph $H$ with vertices $v_1,u_1, u_2, v_2$ on its outer facial cycle in the clockwise order by adding four vertices $s_1, t_1, s_2, t_2$ and ten edges $s_1t_1, s_1u_1, 
s_1v_1, t_1u_1, t_1v_1, s_2t_2, s_2u_2, s_2v_2, t_2u_2, t_2v_2$, where $u_i\ne v_i$ for $i=1,2$ and $|\{u_1,v_1, u_2,v_2\}|\ge 3$ (see Figure \ref{fig:graph in G3}; possibly $\{u_1,v_1\} 
\cap \{u_2,v_2\} \ne \emptyset$).   
\end{itemize}

\vskip 0.5mm

\begin{theorem}[Global structures of Seymour graphs] \label{Thm3} Let $G=(V,E)$ be a graph with two specified nonadjacent edges $e_1=s_1t_1$ and 
$e_2=s_2t_2$. Suppose that $G$ is 2-connected and every 2-separation of $G$, if any, separates $e_1$ from $e_2$. If $G$ is $K_4^*$-free, then  
one of the following statements holds:
\begin{itemize}
\vspace{-2mm}
\item [(i)] $G$ arises from a plane graph $H$ in $\mathcal{G}_1$ by 3-summing 3-connected graphs to its facial triangles.
\vspace{-2mm}
\item [(ii)] $G$ arises from a graph in $\mathcal{G}_2$ by 3-summing 3-connected graphs to triangles $s_1t_1u_1$, $s_1t_1v_1$, and to facial triangles of $H$,
up to interchanging $e_1$ and $e_2$.
\vspace{-2mm}
\item [(iii)] $G$ arises from a graph in $\mathcal{G}_3$ by 3-summing 3-connected graphs to triangles $s_1t_1u_1$, $s_1t_1v_1$, $s_2t_2u_2$, $s_2t_2v_2$, 
and to facial triangles of $H$.
\vspace{-2mm}
\item [(iv)] For each $M$-bridge $B$ of $G$, where $M$ is the graph formed by $e_1$ and $e_2$ only, either $B$ is trivial, or $B$ has exactly three feet, or $B$ has 
a 1-separation $(B_1,B_2)$ with $V(B_1)\cap V(B_2)=\{x\}$ such that $x\notin V(M)$ and $\{s_i,t_i\} \subseteq V(B_i)$ for $i=1,2$.
\vspace{-2mm}
\end{itemize}
\noindent Moreover, the structure described above can be exhibited in $O(n^3)$ time, where $n=|V|$. 
\end{theorem}

The reverse direction of this theorem also holds: If a graph is as described in one of (i)-(iv), then it contains no $K_4^*$-minor either; the proof is 
omitted here, as this statement will not be used in the subsequent sections. 

It is worthwhile pointing out that Truemper \cite{Tr} managed to give a structural description of all matroids with the max-flow min-cut property 
due to Seymour \cite{sey77}, and Gerards et al. \cite{G} obtained a constructive characterization of the so-called odd-$K_4$-free signed graphs. As $K_4^*$-free graphs 
we consider correspond to subclasses of such matroids and such signed graphs, we could in principle deduce Theorem \ref{Thm3} from these existing results.
However, since our graph class is not closed under some splitting operations employed to describe these matroids and signed graphs, such a proof is not so
natural and would require tedious case analysis. Furthermore, in algorithm design efficiency and complexity are always our first concern. Thus we choose to 
present a direct proof without using these existing results.   

Let $P$ be a path (directed or undirected) from $a$ to $b$ and let $c$ and $d$ be two vertices on $P$ such that $a,c,d,b$ (not necessarily distinct) 
occur on $P$ in order as we traverse $P$ in its direction from $a$. We use $P[c,d]$ to denote the subpath of $P$ from $c$ to $d$, and set  
$P(c,d]=P[c,d]\de c$, $P[c,d)=P[c,d]\de d$, and  $P(c,d)=P[c,d]\de \{c,d\}$.

We break the proof of Theorem \ref{Thm3} into a series of lemmas. Throughout $K_4^*$ in respect of a graph $G$ with two specified nonadjacent edges 
$f_1$ and $f_2$ is obtained from one in Figure 2 by replacing $\{e_1,e_2\}$ with $\{f_1, f_2\}$.

\begin{lemma}\label{withtriangle}
Let $G=(V,E)$ be a graph with two specified nonadjacent edges $e_1=s_1t_1$ and $e_2=s_2t_2$. Suppose $G$ has a 3-separation $(G_1,G_2)$ with separating 
set $\{x,y,z\}$, such that $G_1^+:=(V_1,E_1 \cup E(T))$ is 3-connected and that $e_1,e_2\notin E_1\backslash E(T)$, where $G_1=(V_1,E_1)$
and $T$ is the triangle $xyz$. Then $G^+:=(V,E\cup E(T))$ has a $K_4^*$-minor if and only if so does $G$.
\end{lemma}


{\bf Proof.} The ``if" part is trivial. Let us establish the ``only if" part. Suppose $G^+$ has a $K_4^*$-minor (see Figure \ref{fig:K4*}). Then it contains a 
$K_4$-subdivision $J$ with four vertices $a_1,a_2,a_3,a_4$ of degree three, such that $e_1$ is contained in the path $P[a_1,a_2)$ and $e_2$ is contained 
in the path $P[a_3, a_4]$, where $P[a_i,a_j]$ denotes the path between $a_i$ and $a_j$ in $J$ that contains neither of the remaining two 
degree-three vertices, for $1\le i\ne j\le 4$.  As $J$ contains no triangle, at least one of the edges on $T$, say $xy$, is outside $J$.
We propose to show that

(1) No cycle in $J$ is fully contained in $G_1^+\de xy$.

Assume the contrary: Some cycle $C$ in $J$ is fully contained in $G_1^+\de xy$. If neither $e_1$ nor $e_2$ is one $C$, then $C$ 
arises from $J$ by deleting all vertices on $P(a_1,a_2)\cup P(a_3,a_4)$ (see Figure \ref{fig:K4*}). Since $a_1,a_2,a_3,a_4$ are all on $C$ and $e_1,e_2\notin E_1\de E(T)$, 
each of $P[a_1,a_2]$ and $P[a_3,a_4]$ contains at least two vertices from triangle $T$ and hence they are not vertex-disjoint, a contradiction. 
So we may assume the at least one (and hence precisely one) of $e_1$ and $e_2$ is contained in $C$, say $xz=e_1$ and so $e_2 \notin E(T)$.  
Without loss of generality, we may assume that $C$ is $P[a_1,a_2]P[a_2,a_4]P[a_4,a_1]$. Since neither $a_2$ nor $a_4$ is an end of $e_1$, we have 
$\{a_2,a_4\}\cap \{x,z\}=\emptyset$. Note that $e_2 \notin E(G_1^+)$ while $P[a_2,a_3]P[a_3,a_4]$ is a path passing through $e_2$ and internally
vertex-disjoint from $C$, which is impossible. Hence (1) holds. 

We may assume that $J$ is not fully contained in $G$, for otherwise $G$ has a $K_4^*$-minor, we are done. So at least one of $xz$ and $zy$ 
belongs to $E(J)\de E(G)$. From (1) and the structure of $J$, we thus deduce that edges in $J\cap G_1^+$ induce a path, denoted by $Q$, whose two 
ends are contained in $V(T)$. Let $v$ be a vertex in $V(G_1^+)\backslash V(T)$. Since $G_1^+$ is 3-connected, there exist three internally vertex-disjoint paths 
$P_x$, $P_y$, and $P_z$ from $v$ to $x$, $y$, and $z$, respectively, in $G_1^+$ (and hence in $G_1$). It is then a routine matter to check that the 
graph obtained from $J$ by deleting all edges on $Q\de \{e_1,e_2\}$ and adding $P_x \cup P_y \cup P_z$ contains a $K_4^*$-minor in $G$. \qed\\

Let $G$ be a graph as described in Theorem \ref{Thm3}(i), (ii) or (iii). Then it arises from a graph $G_0$ in $\mathcal{G}_1$, $\mathcal{G}_2$, 
or $\mathcal{G}_3$ by 3-summing 3-connected graphs to some triangles in $G_0$; let $\mathcal{T}$ be the collection of all such triangles. We call the
graph $G\cup(\cup_{T\in\mathcal{T}} \,T)$ the {\em closure} of $G$, denoted by $cl(G)$. Repeated applications of Lemma \ref{withtriangle} yield the following
statement.

\begin{corollary}\label{closure}
Let $G$ be a graph as described in Theorem \ref{Thm3}(i), (ii) or (iii). Then $cl(G)$ has a $K_4^*$-minor if and only if so does $G$.
\end{corollary}

We prove Theorem \ref{Thm3} by contradiction. Consider a counterexample $(G,e_1,e_2)$ such that $G$ has the minimum number of edges. Let $M$ be the graph formed 
by $e_1$ and $e_2$ only. By Menger's theorem, we can divide the $M$-bridges of $G$ into four families:
\begin{itemize}
\vspace{-2mm}
\item $\mathcal{B}_2$ consists of all $M$-bridges with two feet; 
\vspace{-2mm}
\item $\mathcal{B}_3$ consists of all $M$-bridges with three feet; 
\vspace{-2mm}
\item $\mathcal{B}_{4a}$ consists of all $M$-bridges $B$ with four feet such that $B$ has a 1-separation $(B_1,B_2)$ separating $\{s_1,t_1\}$ from 
$\{s_2,t_2\}$; and
\vspace{-2mm}
\item $\mathcal{B}_{4b}$ consists of all $M$-bridges $B$ with four feet such that $B$ has two vertex-disjoint paths between $\{s_1,t_1\}$ and $\{s_2,t_2\}$. 
\end{itemize}

\vskip 2mm
\begin{lemma}\label{bridge} The following properties are satisfied by the $M$-bridges of $G$:  
\begin{itemize}
\vspace{-2mm}
\item [(i)] Each $B\in\mathcal{B}_2$ is an edge between $\{s_1,t_1\}$ and $\{s_2,t_2\}$.
\vspace{-2mm}
\item [(ii)] For each $B\in\mathcal{B}_{4a}$, let $(B_1,B_2)$ be a 1-separation of $B$ separating $\{s_1,t_1\}$ from $\{s_2,t_2\}$, with 
$V(B_1)\cap V(B_2)=\{x\}$. Then $x \notin V(M)$.
\vspace{-2mm}
\item [(iii)] $\mathcal{B}_{4b}$ is nonempty.
\vspace{-2mm}
\item [(iv)] Each $B\in\mathcal{B}_{4b}$ contains two vertex-disjoint paths $Q_1$ and $Q_2$ between $\{s_1,t_1\}$ and $\{s_2,t_2\}$,
and also a path $R$ between the internal vertices of $Q_1$ and $Q_2$.
\end{itemize}
\end{lemma}

\vspace{-1mm}
{\bf Proof.} (i) Assume the contrary: $|V(B)|>2$. Let $C$ be the subgraph of $G$ that consists of all vertices outside $B\de V(M)$ 
and all edges outside $B$. Then $(B,C)$ is a 2-separation of $G$ that does not separate $e_1$ from $e_2$, contradicting the hypothesis 
on $G$.

(ii) Otherwise, $B$ would be the union of at least two $M$-bridges of $G$ rather than an $M$-bridge, a contradiction. 

(iii) Suppose not. Then Theorem \ref{Thm3}(iv) would hold true, contradicting the assumption that $(G,e_1,e_2)$ is a counterexample.
     
(iv) By the definition of $\mathcal{B}_{4b}$, bridge $B$ contains two vertex-disjoint paths $Q_1$ and $Q_2$ between $\{s_1,t_1\}$ and $\{s_2,t_2\}$.  
Let $B$ be induced by the edges of a component $C$ of $G\de V(M)$ together with the edges linking $C$ to $V(M)$. Then we may take any path 
in $C$ connecting $Q_1$ and $Q_2$ to be $R$.  \qed

\vskip 4mm

\begin{lemma}\label{B4b}
The family $\mathcal{B}_{4b}$ contains exactly one member $B$. Furthermore, $E(G)\de E(B)$ consists of $e_1$, $e_2$, and $0, 1$, or $2$ edges from $\{f_1,f_2\}$, 
where $f_i$ is an edge between the two ends of $Q_i$ (see Lemma \ref{bridge}(iv)) for $i=1, 2$. 
\end{lemma}

\vspace{-1mm}

{\bf Proof.} By Lemma \ref{bridge}(iii), $\mathcal{B}_{4b}$ is nonempty. Assume on the contrary that there are at least two bridges $B$ and $B'$ in $\mathcal{B}_{4b}$. 
By Lemma \ref{bridge}(iv), $B$ (resp. $B'$) contains two vertex-disjoint paths $Q_1$ and $Q_2$ (resp. $Q_1'$ and $Q_2'$) between $\{s_1, t_1\}$ and $\{s_2, t_2\}$ 
and a path $R$ (resp. $R'$) between the internal vertices of them. Without loss of generality, we may assume that $Q_1$ is an $s_1$-$s_2$ path and $Q_2$ is a 
$t_1$-$t_2$ path. Let $P$ be an $s_1$-$t_2$ path contained in $R' \cup Q_1' \cup Q_2'$. Clearly, $P(s_1,t_2)$ is vertex-disjoint from $B$. So the subgraph of $G$ 
formed by $Q_1$, $Q_2$, $R$, $P$, and $M$ is contractible to $K_4^*$ (see Figure \ref{fig:K4*}); this contradiction implies that $\mathcal{B}_{4b}$ contains exactly one member $B$. 

Similarly, any other $M$-bridge of $G$ consists of only one edge $f_i$ for $i=1$ or $2$.  \qed\\

Let $B$, $Q_1$, $Q_2$, and $R$ be as described in Lemma \ref{B4b} and Lemma \ref{bridge}(iv), respectively, and let $B_0$ be the union of $B$ and $M$. 
Observe that $B_0$ is 2-connected and every 2-separation of $B_0$, if any, separates $e_1$ from $e_2$, because $G$ satisfies these properties. 
Without loss of generality, we assume that $Q_1$ is an $s_1$-$s_2$ path and $Q_2$ is a $t_1$-$t_2$ path (otherwise, we may interchange $s_2$ and $t_2$).

\begin{lemma}\label{1bridge}
$G=B_0$. So $B$ is the only $M$-bridge of $G$.
\end{lemma}

\vspace{-1mm}

{\bf Proof}.  Suppose not. Then Theorem \ref{Thm3} holds for $(B_0,e_1, e_2)$, as $(G, e_1,e_2)$ is a minimum counterexample. Since $B_0$ is a subgraph 
of $G$, it is also $K_4^*$-free. Since $B$ is an $M$-bridge of $B_0$, we see that $B_0$ does not satisfy Theorem \ref{Thm3}(iv). 
    
If $B_0$ satisfies Theorem \ref{Thm3}(i),  then it arises from a plane graph $H$ in $\mathcal{G}_1$ by 3-summing 3-connected graphs to its facial triangles, 
with both $e_1$ and $e_2$ (and hence $s_1,t_1,s_2,t_2$) on the outer facial cycle of $H$. By Lemma \ref{bridge}(iv), the $s_1$-$s_2$ path $Q_1$ and $t_1$-$t_2$ path 
$Q_2$ are vertex-disjoint, so the graph obtained from $H$ by adding $f_1=s_1s_2$ and $f_2=t_1t_2$ remains to be planar. Hence $G$ also satisfies Theorem \ref{Thm3}(i) 
by Lemma \ref{B4b}, a contradiction. 

If $B_0$ satisfies Theorem \ref{Thm3}(ii), then symmetry allows us to assume that it arises from a graph in $\mathcal{G}_2$ by 3-summing 3-connected graphs to 
triangles $s_1t_1u_1$, $s_1t_1v_1$, and to facial triangles of $H$, with $e_2$ and $u_1,v_1$ on the outer facial cycle $C$ of $H$ (see Figure \ref{fig:graph in G2}). Thus the 
subgraph of $cl(B_0)$ (the closure of $B_0$) formed by $\{e_1,e_2\}\cup C\cup\{s_1v_1,t_1u_1\}$ and any of $f_1=s_1s_2$ and $f_2=t_1t_2$ (see Lemma \ref{B4b}) 
is contractible to $K_4^*$. By Corollary \ref{closure}, $B_0$ also has a $K_4^*$-minor, a contradiction. 

If $B_0$ satisfies Theorem \ref{Thm3}(iii), then it arises from a graph in $\mathcal{G}_3$ by 3-summing 3-connected graphs to 
triangles $s_1t_1u_1$, $s_1t_1v_1$, $s_2t_2u_2$, $s_2t_2v_2$ and to facial triangles of $H$, with $u_1,v_1,v_2,u_2$ on the outer 
facial cycle $C$ of $H$ when $H \ne K_2$ or $u_1=u_2$ and $v_1=v_2$ when $H=K_2$. (As usual, $K_n$ stands for the complete graph with $n$ vertices.)
Let $C'$ be the $u_2$-$v_2$ path on $C$ passing $u_1$ and $v_1$ if $H\neq K_2$ or the edge $u_1v_1$ if $H=K_2$. Then the subgraph of $cl(B_0)$ formed 
by $\{e_1,e_2\}\cup C'\cup\{s_1v_1,t_1u_1, s_2u_2,t_2v_2\}$ and any of $f_1=s_1s_2$ and $f_2=t_1t_2$ (see Lemma \ref{B4b}) is contractible to $K_4^*$. 
By Corollary \ref{closure}, $B_0$ also has a $K_4^*$-minor, a contradiction again.   \qed

\vskip 4mm

\begin{lemma}\label{3-connected}
$G$ is 3-connected.
\end{lemma}

\vspace{-1mm}

{\bf Proof.} Assume on the contrary that $G$ has a 2-separation $(G_1,G_2)$. Then it separates $e_1$ from $e_2$ by the hypothesis of the present theorem. Clearly, 
we may assume that $e_i\in E_i$ for $i=1,2$, where $G_i=(V_i,E_i)$. Let $s_3$ and $t_3$ be the two vertices in $V_1\cap V_2$, where $s_3$ is on the $s_1$-$s_2$ 
path $Q_1$ and $t_3$ is on the $t_1$-$t_2$ path $Q_2$. For $i=1,2$, let $G_i^+:=(V_i,E_i \cup\{e_3\})$ and view $e_i$ and $e_3$ as two specified edges in $G_i^+$. 
Observe that

(1) $G_i^+$ is a 2-connected minor of $G$. Furthermore, if $e_i$ is nonadjacent to $e_3$, then $G_i^+$ is $K_4^*$-free, and every 2-separation of $G_i^+$ 
separates $e_i$ from $e_3$. So Theorem \ref{Thm3} holds for $(G_i^+,e_i,e_3)$.

To justify this, note that $G_i^+$ arises from $G_i \cup Q_1 \cup Q_2 \cup \{e_{3-i}\}$ by contracting edges, so it is a 2-connected minor of $G$.
When $e_i$ and $e_3$ are disjoint, $G_i^+$ is $K_4^*$-free, for otherwise, a $K_4^*$-minor in $G_i^+$ can be transformed into a $K_4^*$-minor in 
$G$ by using the $s_3$-$t_3$ subpath of $Q_1 \cup Q_2 \cup \{e_{3-i}\}$. Besides, every 2-separation of $G_i^+$ separates $e_i$ from $e_3$, 
because every 2-separation of $G$ separates $e_1$ from $e_2$.  From the minimality assumption on the counterexample $(G,e_1,e_2)$, we see that 
Theorem \ref{Thm3} holds for $(G_i^+,e_i,e_3)$. Hence (1) is true.

By Lemma \ref{1bridge}, $B$ is the only $M$-bridge of $G$, so at most one of $s_3$ and $t_3$ comes from $V(M)$. For convenience, we choose a 2-separation 
$(G_1,G_2)$ of $G$ so that 

(2) $\{s_3,t_3\} \cap V(M)\ne \emptyset$ if possible (now we may assume that $s_3=s_1$) and $V_1$ is inclusionwise minimal otherwise. 

When $s_3=s_1$, let $G_1^*$ be the union of $G_1^+$ and edge $t_1t_3$. It is a routine matter to check that $G_1^*$ is $3$-connected, because every 
2-separation of $G$ separates $e_1$ from $e_2$. Thus 

(3) If $s_3=s_1$, then $G_1^+$ is the 3-sum of a 3-connected graph $G_1^*$ and the triangle $s_1t_1t_3$.

(4) If $s_3\neq s_1$, then $G_1^+$ is as described in Theorem \ref{Thm3}(i) or (iv).

To justify this, observe that if $G_1^+$ is as described in Theorem \ref{Thm3}(ii) or (iii), then $\{u_1,v_1\}$ or $\{u_2,v_2\}$ (see Figures \ref{fig:graph in G2} and \ref{fig:graph in G3})
would be a $2$-separating set of $G$ that separates $e_1$ from $e_2$, contradicting the minimality assumption on $V_1$ (see (2)). So (4) holds.

(5) $G_2^+$ is as described in Theorem \ref{Thm3}(i), or (ii) with $e_3$ on the outer facial cycle $C_2$ of $H$ (see Figure \ref{fig:graph in G2}), or (iv).

To prove this, recall that if $G_2^+$ is as described in Theorem \ref{Thm3}(ii), then it arises from a graph $H'_2$ in $\mathcal{G}_2$ by 3-summing 3-connected 
graphs to some triangles of $H'_2$, up to interchanging $e_2$ and $e_3$ (see the definition of $\mathcal{G}_2$). So either $e_2$ or $e_3$ 
(but not both) is on the outer facial cycle $C_2$ of $H_2$. Assume on the contrary that $e_2$, $u_1$ and $v_1$ are on $C_2$. Renaming $u_1$ and $v_1$ if necessary, 
we may assume that $u_1$ is on $Q_1$ and $v_1$ is on $Q_2$. Then the graph formed by $\{e_1,e_2\}\cup Q_1\cup Q_2\cup\{s_3v_1,t_3u_1\}$ 
is contractible to $K_4^*$. Repeated applications of Lemma \ref{withtriangle} imply that $G$ also has a $K_4^*$-minor, a contradiction. 

Similarly, if $G_2^+$ is as described in Theorem \ref{Thm3}(iii), then $G$ also contains $K_4^*$ as a minor; this contradiction establishes (5).

Observe that if $G_1^+$ is as described in (3) or in Theorem \ref{Thm3}(i) and $G_2^+$ is as described in Theorem \ref{Thm3}(i) or (ii), with $e_3$ 
is on $C_2$, then clearly $G$ satisfies Theorem \ref{Thm3}(i) or (ii), contradicting the assumption that $(G,e_1,e_2)$ is a counterexample to
Theorem \ref{Thm3}. Thus (3)-(5) and symmetry allow us to assume hereafter that one of the 
following four cases occurs:

\vskip 1mm
{\bf Case 1.} $G_1^+$ is as described in (3) and  $G_2^+$ is as described in Theorem \ref{Thm3}(iv);

{\bf Case 2.} $G_1^+$ is as described in Theorem \ref{Thm3}(i) and $G_2^+$ is as described in Theorem \ref{Thm3}(iv);

{\bf Case 3.} $G_1^+$ is as described in Theorem \ref{Thm3}(iv) and $G_2^+$ is as described in Theorem \ref{Thm3}(ii), with $e_3$ on $C_2$.

{\bf Case 4.} $G_1^+$ and $G_2^+$ are both as described in Theorem \ref{Thm3}(iv).

\vskip 1mm

Assume that $G_i^+$ is as described in Theorem \ref{Thm3}(iv) for $i=1$ or $2$. Let $M_i$ denote the graph formed by $e_3$ and $e_i$ only. Then the $M_i$-bridges 
of $G_i^+$ can be divided into three families:
\begin{itemize}
\vspace{-2mm}
\item $\mathcal{B}_2^i$ consists of all $M_i$-bridges with two feet; 
\vspace{-2mm}
\item $\mathcal{B}_3^i$ consists of all $M_i$-bridges with three feet; and
\vspace{-2mm}
\item $\mathcal{B}_{4a}^i$ consists of all $M_i$-bridges $J$ with four feet, such that $J$ has a 1-separation $(J_1,J_2)$ separating  $\{s_3,t_3\}$ from 
$\{s_i,t_i\}$, with $V(J_1)\cap V(J_2)=\{g\}$ and $g \notin V(M_i)$.
\end{itemize}
\vspace{-2mm}

(6) For any three feet $a_1,a_2,a_3$ of a bridge $J \in \mathcal{B}_3^i \cup \mathcal{B}_{4a}^i$, there exist three internally vertex-disjoint paths
$P_{a_1}, P_{a_2}, P_{a_3}$ in $J$ from some vertex $b \in V(J)\de V(M_i)$ to $a_1,a_2,a_3$, respectively, such that none of these paths contains
the vertex in $V(M_i)\de \{a_1,a_2,a_3\}$. 

Let $C$ be a component of $G_i^+\de V(M_i)$, such that $J$ is a subgraph of $G_i^+$ induced by the edges of a component $C$ of $G_i^+\de V(M_i)$ 
together with the edges linking $C$ to $M_i$. Let $c_1,c_2,c_3$ (not necessarily distinct) be neighbors of $a_1,a_2,a_3$ in $C$, respectively, and let $T$ be
an inclusivewise minimal tree containing $c_1,c_2,c_3$ in $C$. Obviously, $T \cup \{a_1c_1,a_2c_2,a_3c_3\}$ contains three internally vertex-disjoint paths
$P_{a_1}, P_{a_2}, P_{a_3}$ with the desired properties. So (6) is justified. 

(7) Assume $J$ is a bridge in $\mathcal{B}_{4a}^i$. Then for each $a \in \{s_i,t_i\}$, there exists $b \in \{s_3,t_3\}$ such that $G_i^+$ contains no
path between $a$ and $b$, with all internal vertices outside $J$. 

Otherwise, there exist paths $L_1$ and $L_2$ in $G_i^+$ from $a$ to $s_3$ and $t_3$, respectively, with internal vertices outside $J$. Clearly, we
may assume that $L_1$ and $L_2$ have only one common subpath. Let $c$ be the vertex in $V(M_i)\de \{s_3,t_3,a\}$. Then (6) guarantees the existence of
three internally vertex-disjoint paths $P_{s_3}, P_{t_3}, P_c$ in $J$ from some vertex $x \in V(J)\de V(M_i)$ to $s_3,t_3,c$, respectively, such that 
none of these paths contains $a$. Then the subgraph formed by $\{e_1,e_2\}\cup Q_1[s_{3-i},s_3]\cup Q_2[t_{3-i},t_3]\cup \{L_1, L_2, P_{s_3}, P_{t_3},
P_c\}$ is contractible to $K_4^*$. By Corollary \ref{closure}, $G$ also has a $K_4^*$-minor; this contradiction proves (7).  

(8) Assume $J$ is a bridge in $\mathcal{B}_{4a}^i$ that is vertex-disjoint from $Q_2(t_3,t_i)$, or assume $\mathcal{B}_{4a}^i=\emptyset$ and 
$J$ is a bridge in $\mathcal{B}_3^i$ with feet $s_3$, $t_3$ and $s_i$. Then $G_i^+$ contains no path between $s_3$ and $t_i$, with all internal 
vertices outside $J$.  

Suppose on the contrary that $G_i^+$ contains a path $L$ between $s_3$ and $t_i$, with internal vertices outside $J$. Clearly, we
may assume that $L$ and $Q_2[t_3,t_i]$ have only one common subpath. Note that if $\mathcal{B}_{4a}^i=\emptyset$, then $J$ in $\mathcal{B}_3^i$ 
is vertex-disjoint from $Q_2(t_3,t_i]$. By (6), there exist three internally vertex-disjoint paths $P_{s_3}, P_{t_3}, P_{s_i}$ in $J$ from 
some vertex $x \in V(J)\de V(M_i)$ to $s_3,t_3, s_i$, respectively, such that none of these paths contains $t_i$. Thus the subgraph formed by 
$\{e_1,e_2\}\cup Q_1[s_{3-i},s_3]\cup Q_2\cup \{L, P_{s_3}, P_{t_3}, P_{s_i}\}$ is contractible to $K_4^*$. By 
Corollary \ref{closure}, $G$ also has a $K_4^*$-minor; this contradiction proves (8).  

The following statement is just a mirror image of (8).

(9) Assume $J$ is a bridge in $\mathcal{B}_{4a}^i$ that is vertex-disjoint from $Q_1(s_3,s_i)$, or assume $\mathcal{B}_{4a}^i=\emptyset$ and 
$J$ is a bridge in $\mathcal{B}_3^i$ with feet $s_3$, $t_3$ and $t_i$. Then $G_i^+$ contains no path between $t_3$ and $s_i$, with all internal 
vertices outside $J$.  

Observe that every bridge in $\mathcal{B}_2^i$ is trivial, for otherwise, $G$ would have a $2$-separation that does not separate $e_1$ from $e_2$,
a contradiction. By (7), there is at most one bridge in $\mathcal{B}_{4a}^i$.  By the definition of $\mathcal{B}_{4a}^i$, the bridge $J$ in $\mathcal{B}_{4a}^i$,
if any, is vertex-disjoint from $Q_2(t_3,t_i)$ or vertex-disjoint from $Q_1(s_3,s_i)$.  From (7), (8) and (9), we conclude that $G_i^+$ may have 
the following four types of structure:

{\bf Type 1.}  $\mathcal{B}_{4a}^i$ contains exactly one bridge $J$, and exactly one of (1.1) and (1.2) holds:  

(1.1) all bridges in $\mathcal{B}_2^i \cup \mathcal{B}_3^i$ have feet in the set 
$\{s_3,s_i,t_i\}$; 

(1.2) all bridges in $\mathcal{B}_2^i \cup \mathcal{B}_3^i$ have feet in the set $\{t_3,s_i,t_i\}$. 

{\bf Type 2.} $\mathcal{B}_{4a}^i$ contains exactly one bridge $J$, $\mathcal{B}_3^i=\emptyset$, and  $\mathcal{B}_2^i \subseteq \{s_3s_i, t_3t_i\}$.

{\bf Type 3.} $\mathcal{B}_{4a}^i=\emptyset$, and exactly one of (3.1) and (3.2) holds:  

(3.1) the feet set of each bridge in $\mathcal{B}_3^i$ is either $\{s_3,t_3,s_i\}$ (this case must occur)
or $\{t_3,s_i,t_i\}$, and $s_3t_i \notin \mathcal{B}_2^i$;

(3.2) the feet set of each bridge in $\mathcal{B}_3^i$ is either $\{s_3,t_3,t_i\}$ (this case must occur) or $\{s_3,s_i,t_i\}$, and $t_3s_i 
\notin \mathcal{B}_2^i$.

{\bf Type 4.} $\mathcal{B}_{4a}^i=\emptyset$, the feet set of each bridge in $\mathcal{B}_3^i$ is either $\{s_3,s_i,t_i\}$ or $\{t_3,s_i,t_i\}$. (There
is no restriction on $\mathcal{B}_2^i$.)

\begin{figure}[ht]
    \centering
        \begin{subfigure}{0.4\textwidth}
        \centering
        \includegraphics{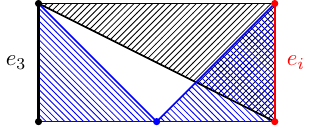} 
        \caption*{$G_i^+$ of Type 1}
        \end{subfigure}
        \begin{subfigure}{0.4\textwidth}
        \centering
        \includegraphics{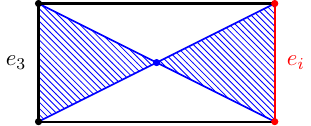}
        \caption*{$G_i^+$ of Type 2}
        \end{subfigure}\\[2ex]
    
        \begin{subfigure}{0.4\textwidth}
        \centering
        \includegraphics{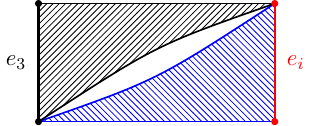}
        \caption*{$G_i^+$ of Type 3}
        \end{subfigure}
        \begin{subfigure}{0.4\textwidth}
        \centering
        \includegraphics{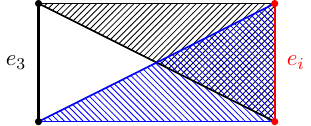}
        \caption*{$G_i^+$ of Type 4}
        \end{subfigure}
    
        \caption{Four types of $G_i^+$}
        \label{fig:T1-T4}
\end{figure}

Figure \ref{fig:T1-T4} illustrates the graph $G_i^+$ of Types 1-4, in which each shaded part with the same color represents a nontrivial $M_i$-bridges. Recall 
the definition of $\mathcal{B}_{4a}^i$; bridge $J$ has a 1-separation $(J_1,J_2)$ separating $\{s_3,t_3\}$ from $\{s_i,t_i\}$, with $V(J_1)\cap V(J_2)=\{g\}$ 
and $g \notin V(M_i)$.

(10) If $G_i^+$ is of Type 1, letting $H_i$ be the graph consisting of three triangles $s_3t_3g$, $s_it_is_3$ and $s_it_ig$ if (1.1) holds and consisting of   
three triangles $s_3t_3g$, $s_it_it_3$ and $s_it_ig$ if (1.2) holds, then $G_i^+$ arises from $H_i$ by 3-summing 3-connected graphs to these three triangles 
in $H_i$. 

(11) If $G_i^+$ is of Type 2, letting $H_i$ be the graph consisting of four triangles $s_3t_3g$, $s_it_ig$, $s_3s_ig$ and $t_3t_ig$, then $G_i^+$ arises from $H_i$ 
by 3-summing 3-connected graphs to facial triangles $s_3t_3g$ and $s_it_ig$ in $H_i$. 

(12) If $G_i^+$ is of Type 3, letting $H_i$ be the graph consisting of two triangles $s_3t_3s_i$ and $s_it_it_3$ if (3.1) holds and consisting of two triangles
$s_3t_3t_i$ and $s_it_is_3$ if (3.2) holds, then $G_i^+$ arises from $H_i$ by 3-summing 3-connected graphs to these two facial triangles in $H_i$. 

(13) If $G_i^+$ is of Type 4, letting $H_i$ be the graph consisting of the edge $e_3$ and two triangles $s_it_is_3$ and $s_it_it_3$, then $G_i^+$ arises from 
$H_i$ by 3-summing 3-connected graphs to these two triangles in $H_i$.

Let us now analyze Cases 1-4 described above one by one.

{\bf Case 1.} $G_1^+$ is as described in (3) and $G_2^+$ is as described in Theorem \ref{Thm3}(iv). In this case, $s_1=s_3$. So 

(14) $G_2^+$ neither contains $s_1s_2$ as a trivial $M_2$-bridge nor contains a nontrivial $M_2$-bridge with feet set $\{s_1,s_2,t_2\}$, 
because $B$ is the only $M$-bridge of $G$ by Lemma \ref{1bridge}.  

(15) $G_2^+$ is not of Type 4.

Otherwise, $\mathcal{B}_{4a}^2=\emptyset$ and the feet set of each bridge in $\mathcal{B}_3^2$ is either $\{s_3,s_2,t_2\}$ or $\{t_3,s_2,t_2\}$.
Let $K$ denote the nontrivial $M_2$-bridge of $G_2^+$ containing $Q_1$ (see (14)). Then $t_3 \notin V(K)$. It follows from (14) that $t_2 \notin V(K)$. Hence
$\{s_1,s_2\}$ would be a $2$-separating set in $G$ that does not separate $e_1$ from $e_2$, contradicting the hypothesis of Theorem \ref{Thm3}. So
(15) is true. 

In view of (15), let $H_2$ be as specified in (10), (11) or (12) and let $H'$ be the union of $H_2$ and the triangle $s_1t_1t_3$. From (3) and (10)-(12), we see that 
$H'\in {\cal G}_1 \cup {\cal G}_2$. Furthermore, $G$ arises from $H'$ by 3-summing 3-connected graphs to certain triangles and satisfies 
Theorem \ref{Thm3}(i) or (ii).

{\bf Case 2.} $G_1^+$ is as described in Theorem \ref{Thm3}(i) and $G_2^+$ is as described in Theorem \ref{Thm3}(iv). In this case, $G_1^+$ arises from a 
plane graph $H_1$ in $\mathcal{G}_1$ by 3-summing 3-connected graphs to its facial triangles. Let $H_2$ be as specified in (10)-(13) and let $H'=H_1 \cup H_2$.
From (10)-(13), we see that $H'\in {\cal G}_1 \cup {\cal G}_2$. Furthermore, $G$ arises from $H'$ by 3-summing 3-connected graphs to certain 
triangles and satisfies Theorem \ref{Thm3}(i) or (ii).

{\bf Case 3.} $G_1^+$ is as described in Theorem \ref{Thm3}(iv) and $G_2^+$ is as described in Theorem \ref{Thm3}(ii), with $e_3$ on $C_2$. In this case,
$G_1$ is not of Type 1, for otherwise, $\{s_1,g\}$ if (1.1) holds and $\{t_1,g\}$ if (1.2) holds would be a $2$-separating set of $G$ that separates $e_1$ 
from $e_2$, contradicting the minimality hypothesis on $|V_1|$ (see (2)). By an argument analogous to that used for the above cases, we can also deduce that 
$G$ satisfies Theorem \ref{Thm3}(ii) or (iii). 
 
{\bf Case 4.} $G_1^+$ and $G_2^+$ are both as described in Theorem \ref{Thm3}(iv). In this case, once again $G_1$ is not of Type 1. From (10)-(13), we can
similarly deduce that $G$ satisfies Theorem \ref{Thm3}(i), (ii) or (iii). 
 
Therefore, in each case we can reach a contradiction to the assumption that $(G,e_1,e_2)$ is a counterexample to
Theorem \ref{Thm3}. This completes the proof of our lemma. \qed\\

Let $Z\subseteq V(G)$. We say that $(G,Z)$ is 4-connected if for every $k$-separation $(G_1,G_2)$ of $G$ with $Z\subseteq V(G_1)$, either $k\geq 4$ or $k=3=|V(G_2)|-1$.

\begin{lemma}\label{4-connected}
$(G,V(M))$ is 4-connected. 
\end{lemma}
\vspace{-1mm}

{\bf Proof.} Assume the contrary: $G$ has a 3-separation $(G_1,G_2)$ with $M$ contained in $G_1$ and $|V(G_2)|>4$. Let us choose such a separation that 
$G_2$ is inclusionwise maximal. Write $V(G_1)\cap V(G_2)=\{x,y,z\}$. We use $G_i^+$ to denote the union of $G_i$ and the triangle $T:=xyz$ for $i=1,2$.
Then $G$ is a 3-sum of $G_1^+$ and $G_2^+$ over $T$ and each $G_i^+$ is 3-connected by Lemma \ref{3-connected}. In view of Lemma \ref{withtriangle},
$G_1^+$ is $K_4^*$-free. So Theorem \ref{Thm3} holds for $(G_1^+, e_1,e_2)$.

If $G_1^+$ satisfies Theorem \ref{Thm3}(i),  then it arises from a plane graph $H$ in $\mathcal{G}_1$ by 3-summing 3-connected graphs to its facial triangles.
Observe that either $T$ is contained in a 3-connected graph $K$ which is 3-summed to some facial triangle of $H$, or $T$ is contained in $H$. In the former
case, $G$ also satisfies Theorem \ref{Thm3}(i), as 3-summing $G_2^+$ to $T$ in $K$ preserves the 3-connectivity of $K$. In the latter case, $T$ must be a facial 
triangle of $H$ by the maximality assumption on $G_2$. Thus $G$ also satisfies Theorem \ref{Thm3}(i) by 3-summing $G_2^+$ to $T$.

If $G_1^+$ satisfies Theorem \ref{Thm3}(ii) or (iii), then it arises from a graph $H_1'$ in $\mathcal{G}_2$ or $\mathcal{G}_3$ by 3-summing 3-connected 
graphs to some triangles in $H$ (see the definitions). In either case, $\{u_1,v_1\}$ or $\{u_2,v_2\}$ (see Figures \ref{fig:graph in G2} and \ref{fig:graph in G3}) would be a 2-separating set of $G$,
contradicting Lemma \ref{3-connected}.

If $G_1^+$ satisfies Theorem \ref{Thm3}(iv), then $V(T)\nsubseteq V(M)$ and $G_1^+$ has exactly one nontrivial $M$-bridge $J$, with feet set $V(M)$,
because $B$ is the only $M$-bridge of $G$ by Lemma \ref{1bridge}. According to Theorem \ref{Thm3}(iv), $J$ has a 1-separation $(J_1,J_2)$ with

(1) $V(J_1)\cap V(J_2)=\{h\}$, $h\notin V(M)$, and $\{s_i,t_i\}\subseteq V(J_i)$ for $i=1,2$. 

\noindent Observe that $J$ can be formed by 3-summing 3-connected graphs to the union of two triangles $s_1t_1h$ and $s_2t_2h$. Furthermore, 
$|V(T)\cap V(M)|=2$, for otherwise, some 3-separation of $G$ with separating set $\{s_i,t_i,h\}$ for $i=1$ or $2$ would violate the maximality 
assumption on $G_2$. If $\{s_i,t_i\}\subset V(T)$ for $i=1$ or $2$, then $T$ is contained in $J_i$ and hence $G$ also satisfies Theorem \ref{Thm3}(iv).
So $V(T)\cap \{s_i,t_i\}\neq\emptyset$ for $i=1,2$. It follows from (1) that the third vertex of $T$ must be $h$. Thus $G$ can be obtained by 3-summing 
3-connected graphs to the union of three triangles $s_1t_1h$, $s_2t_2h$ and $T$, and hence satisfies Theorem \ref{Thm3}(i).
 
Therefore, we can reach a contradiction in any case, which implies that $(G,V(M))$ is 4-connected.   \qed\\

In our proof we also need a theorem from Ding and Marshall \cite{DM}, which can be used to decide if a given graph admits a planar drawing with some distinguished 
vertices and edges lying on the same facial cycle. A {\it circlet} $\Omega$ of a graph $G$ consists of a cyclically ordered set of distinct vertices $v_1,v_2,\ldots,v_n$ 
of $G$, where $n\geq 4$, and a set of edges of $G$ of the form $v_iv_{i+1}$, where $v_{n+1}=v_1$. Note that not necessarily all edges of $G$ of the given form are in $\Omega$. 
An {\em $\Omega$-cycle} is a cycle $C$ of $G$ such that $V(\Omega)\subseteq V(C)$, $E(\Omega)\subseteq E(C)$, and the cyclic ordering of $V(\Omega)$ 
agrees with the ordering in $C$. For each $i$, the $v_i$-$v_{i+1}$ path of $C$ that does not contain $v_{i+2}$ is called a {\em segment} of $C$. 
For a cycle $C$ in $G$, a path $P$ of $G$ is called a {\em $C$-path} if $E(P)\cap E(C)=\emptyset$ and the distinct ends of $P$ are the only two 
vertices of $P$ that are in $C$. (So no edge on $C$ is a $C$-path.) A pair of vertex-disjoint $C$-paths $(P_1, P_2)$ is called a {\em cross} of $C$ 
if the ends of $P_1$, say $x_1$ and $y_1$, and the ends of $P_2$, say $x_2$ and $y_2$, appear in the order $x_1, x_2, y_1, y_2$ around $C$.

\begin{theorem}[Theorem 7.3 in Ding and Marshall \cite{DM}]\label{planar drawing}
Let $\Omega$ be a circlet of $G$ such that $G$ has an $\Omega$-cycle and $(G, V(\Omega))$ is 4-connected. Then either $G$ admits a planar drawing in which some 
facial cycle is an $\Omega$-cycle, or $G$ has an $\Omega$-cycle $C$ and a cross of $C$ for which each segment of $C$ contains at most two of the four ends of 
the two paths in the cross.
\end{theorem}

Now we are ready to finish our structural description of Seymour graphs.
\vskip 2mm

{\bf Proof of Theorem \ref{Thm3}.} Lemma \ref{4-connected} allows us to apply Theorem \ref{planar drawing} for the circlet $\Omega$ with cyclically ordered 
vertex set $\{s_1,s_2,t_2,t_1\}$ and edge set $\{e_1,e_2\}$, as the cycle formed by $Q_1$, $Q_2$, $e_1$, and $e_2$ is an $\Omega$-cycle. If $G$ 
admits a planar drawing in which some facial cycle is an $\Omega$-cycle, then $G$ is as described in Theorem \ref{Thm3}(i), contradicting the assumption 
that $(G_1, e_1,e_2)$ is a counterexample. Thus, by Theorem \ref{planar drawing}, $G$ has an $\Omega$-cycle $C$ and a cross $(P_1,P_2)$ of $C$ 
for which each segment of $C$ contains at most two of the four ends of $P_1$ and $P_2$. Let $x_i$, $y_i$ be the ends of $P_i$ for $i=1,2$. Clearly,
we may assume that $s_1,x_1,x_2, s_2,t_2, y_1, y_2, t_1$ occur on $C$ in the clockwise order. For any two vertices $a, b$ on $C$, let $C[a,b]$ be the 
path of $C$ from $a$ to $b$ in the clockwise direction.  We use $H$ to denote the subgraph of $G$ comprising $C$, $P_1$, and $P_2$.  

By Lemma \ref{1bridge}, $B$ is the only $M$-bridge of $G$. So

(1) $s_1s_2, t_1t_2 \notin E(G)$.

Note that if neither $\{s_1,t_1\}=\{x_1,y_2\}$ nor $\{s_2,t_2\}=\{x_2,y_1\}$ holds true, then the subgraph of $G$ formed by $C\cup P_1\cup P_2$ is 
contractible to $K_4^*$, a contradiction. So 

(2) $\{s_1,t_1\}=\{x_1,y_2\}$ or $\{s_2,t_2\}=\{x_2,y_1\}$.

\noindent We may assume that  

(3) $\{x_1,y_1,x_2,y_2\} \ne V(M)$. 

Otherwise,  $x_1=s_1$, $y_2=t_1$, $x_2=s_2$, and $y_1=t_2$. By Lemma \ref{3-connected}, $G$ is 3-connected, so 
$G\de \{s_1,s_2\}$ contains a path $P_3$ connecting a vertex $x_3$ on $C(s_1,s_2)$ (see (1)) and a vertex $y_3$ in $H\de V(C[s_1,s_2])$, such that 
the internal vertices of $P_3$ are all outside $H$. If $P_3$ is vertex-disjoint from $P_1 \cup P_2$, then $y_3$ is on $C(t_2,t_1)$. Thus the 
subgraph $C\cup P_1\cup P_3$ of $G$ is contractible to $K_4^*$, a contradiction. So $P_3$ intersects $P_1$ or $P_2$, say $y_3\in V(P_1)$. 
Then $P_3[x_3,y_3]P_1[y_3,y_1]$ and $P_2$ form another cross of $C$ whose four ends are not all contained in $V(M)$. Hence we may assume (3).

Symmetry, (2) and (3) allow us to assume that the $\Omega$-cycle $C$ and its cross $(P_1,P_2)$ are selected so that 

(4) $\{s_1,t_1\}=\{x_1,y_2\}$ and $\{s_2,t_2\}\ne \{x_2,y_1\}$;  

(5) subject to (4), the path $C[x_2,y_1]$ is as short as possible. 

Depending on the locations of $x_2$ and $y_1$ on $C$, we distinguish between two cases. 

{\bf Case 1.} $\{s_2,t_2\}\cap\{x_2,y_1\}\neq\emptyset$, say $x_2=s_2$ (so $y_1\ne t_2$ by (4)).  In this case, let $P_3$ be a path connecting 
a vertex $x_3$ on $C(s_1,s_2)$ and a vertex $y_3$ in $H\de V(C[t_1,t_2])$, such that the internal vertices of $P_3$ are all outside $H$; this
path $P_3$ is available because, by Lemma \ref{1bridge}, $B$ is the only $M$-bridge in $G$ and hence $G\backslash V(M)$ is connected. Now the same 
argument used in the proof of (3) implies that $H \cup P_3$ contains a $K_4^*$-minor, a contradiction.  

{\bf Case 2.} $\{s_2,t_2\}\cap\{x_2,y_1\}=\emptyset$. In this case, let $P_3$ be a path in $G\backslash\{x_2,y_1\}$ connecting a vertex $x_3$ on $C(x_2,y_1)$
and a vertex $y_3$ in $H \backslash V(C[x_2,y_1])$, such that the internal vertices of $P_3$ are all outside $H$; this path $P_3$ is available because $G$ 
is 3-connected. By symmetry, we may assume that $x_3$ is on $C(x_2,s_2]$. There are four possibilities for the location of $y_3$: 

$\bullet$ $y_3$ is on $C[s_1, x_2)$. In this subcase, let $C'$ be obtained from $C$ by replacing $C[y_3,x_3]$ with $P_3$ and let $P_2'$ the path 
$P_2 \cup C[x_2,x_3]$. Then the existence of the $\Omega$-cycle $C'$ and its cross $(P_1,P_2')$ would contradict the choice of $\Omega$-cycle 
$C$ and its cross $(P_1,P_2)$, as (5) is violated. 

$\bullet$ $y_3$ is on $C(y_1,t_1]$. In this subcase,  if $y_3=t_1$, then the existence of the $\Omega$-cycle $C$ and its cross $(P_1,P_3)$ would 
contradict the choice of $\Omega$-cycle $C$ and its cross $(P_1,P_2)$, as (5) is violated. If $y_3\ne t_1$, then the subgraph $C \cup P_1 \cup P_3$
of $G$ would be contractible to $K_4^*$, a contradiction.   

$\bullet$ $y_3$ is on $P_1(x_1,y_1)$. In this subcase,  $H\de V(C(s_1,x_2))$ would be contractible to $K_4^*$, a contradiction.   

$\bullet$ $y_3$ is on $P_2(x_2,y_2)$. In this subcase, let $P_2'=P_3\cup P_2[y_3,y_2]$. Then the existence of the $\Omega$-cycle $C$ and 
its cross $(P_1,P_2')$ would contradict the choice of $\Omega$-cycle $C$ and its cross $(P_1,P_2)$, as (5) is violated.

Combining the above cases, we conclude that $(G,e_1,e_2)$ cannot be a counterexample to Theorem \ref{Thm3}. So $G$ is indeed as described in
(i)-(iv). 

Although the whole proof proceeds by reductio ad absurdum, it can be easily converted to an $O(n^3+m)$ algorithm for finding the global structure
of $G$, where $n=|V(G)|$ and $m=|E(G)|$. Let us first find all connected components of $G\de V(M)$ by depth-first search, thereby obtaining all nontrivial 
$M$-bridges of $G$ accordingly; this step takes $O(n+m)$ time (see Hopcroft and Tarjan \cite{HT3}). If no nontrivial bridge with four feet has two vertex-disjoint 
paths between $\{s_1,t_1\}$ and $\{s_2,t_2\}$, then $G$ is as described in Theorem \ref{Thm3}(iv). So we assume that $B$ is such a bridge. Let $B_0$ 
be the union of $G$ and $M$. By Lemma \ref{B4b},  $G$ is obtained from $B_0$  by adding 0, 1, or 2 edges from $\{f_1,f_2\}$. By induction, the structure 
of $B_0$ can be found in time $O(n^3+|E(B_0)|)$ time and hence that of $G$ can be obtained 
in $O(n^3+m)$ time (see the proof of Lemma \ref{1bridge}). So we assume that $G=B_0$. If $G$ is not $3$-connected, then we can find all $3$-connected
components and $2$-separating sets again by depth-first search, which requires $O(n+m)$ time (see Hopcroft and Tarjan \cite{HT1}). Consequently, we can determine
a $2$-separation $(G_1,G_2)$ and hence $G_1^+$ and $G_2^+$ as specified in the proof of Lemma \ref{3-connected}. By induction, the structure of $G_i^+$ 
can be obtained in time $O(n_i^3+m_i)$ time, where $n_i$ and $m_i$ are the numbers of vertices and edges in $G_i^+$ respectively. Thus the structure
of $G$ can be exhibited in $O(n+m)+O(n_1^3+m_1)+O(n_2^3+m_2)=O(n^3+m)$ time. So we further assume that $G$ is $3$-connected.  Using Kanevsky and 
Ramachandran's algorithm \cite{KR}, we can output all $3$-separating sets of $G$ in $O(n^2)$ time, thereby determining if $(G, V(M))$ is 
$4$-connected. If not, let $(G_1, G_2)$ be the $3$-separation as specified in the proof of Lemma \ref{4-connected} and find the structure of $G_1^+$
is $O(n_1^3+m_1)$ time, where $n_i$ and $m_i$ are the numbers of vertices and edges in $G_i^+$ respectively for $i=1,2$. Thus the structure of $G$
can be found in $O(n^2)+O(n_1^3+m_1)+O(n_2+m_2)=O(n^3+m)$ time. Finally, assume $(G, V(M))$ is $4$-connected. Then $G$ admits a planar drawing in which 
both $e_1$ and $e_2$ are on its outer facial cycle, so it is as described in Theorem \ref{Thm3}(i) (see the beginning of this proof); such a drawing can 
be obtained in $O(n+m)$ time using Hopcroft and Tarjan's algorithm \cite{HT2}, completing the proof. \qed

\section{Network Flows} 

In Section 1, we have introduced the biflows in the {\it path packing} form, which is concise in language and more convenient to use from the
combinatorial viewpoint. However, when it comes to algorithm design, we often need to represent flows and biflows in the {\it arc-vertex} form. 

Let $G=(V,E)$ be a graph with a capacity function $c:E \rightarrow \mathbb{Z}_+$. We use $D=(V,A)$ to denote the digraph obtained from $G$ by
replacing each edge $uv \in E$ with a pair of opposite arcs $(u,v)$ and $(v,u)$, such that $c(u,v)=c(v,u)=c(uv)$. For arc $a=(u,v)$, its {\em reverse},
denoted by $\bar{a}$, is $(v,u)$. For each $U \subseteq V$, let $\delta(U)$ denote the set of all edges of $G$ between $U$ and $V\de U$, let 
$\delta^+(U)$ denote the set of all arcs of $D$ from $U$ to $V\de U$, and let $\delta^-(U)$ denote the set of all arcs of $D$ from $V\de U$ to $U$. 
When $U=\{u\}$, we write $\delta(u)$ for $\delta(U)$, $\delta^+(u)$ for $\delta^+(U)$, and $\delta^-(u)$ for $\delta^-(U)$. 

Let $S$ and $T$ be two nonempty disjoint subsets of $V$. An $(S,T)$-{\em cut} in the network $(G,c)$ (or simply $G$) is a subset $K$ of $E$ 
that contains at least one edge from each path connecting $S$ and $T$, whose capacity is defined to be $c(K)$. Clearly, every inclusionwise minimal 
$(S,T)$-cut is of the form $\delta(U)$ for some vertex subset $U$ with $S \subseteq U$ and $T \subseteq V\de U$. So to find a minimum $(S,T)$-cut 
(with respect to capacity), we may restrict our attention to such cuts. An $(S,T)$-{\em flow} in the network $(G,c)$ (or simply $G$) is 
a function $f: A \rightarrow \mathbb{R}_+$ such that 

$(3.1)$ $0\le f(a) \le c(a)$ for each arc $a \in A$;

$(3.2)$ $f(a)+f(\bar{a}) \le c(a)$ for each arc $a \in A$; and

$(3.3)$ ${\rm div}_f(v)=0$ if $v \in V\de (S\cup T)$, ${\rm div}_f(v)\ge 0$ if $v \in S$, and ${\rm div}_f(v)\le 0$ if $v \in T$, 

\noindent where ${\rm div}_f(v) = f(\delta^+(v))- f(\delta^-(v))$. It is easy to see that ${\rm div}_f(S)=-{\rm div}_f(T)$; this quantity is called 
the {\em value} of $f$, denoted by ${\rm val}(f)$. When $S=\{s\}$ (resp. $T=\{t\}$), we write $(s,T)$ (resp. $(S,t)$) for $(S,T)$. When $S=\{s\}$ 
and $T=\{t\}$, an $(S,T)$-cut is often called an $s$-$t$ {\em cut}, and an $(S,T)$-flow is often called an $s$-$t$ {\em flow}. Note that 
a flow $f$ in the undirected network $(G,c)$ is actually a flow in the directed network $(D,c)$, and that in the classical network flow problems 
on directed networks, there is no constraint $(3.2)$. Furthermore, we assume that 
 
$(3.4)$ $f(a)f(\bar{a})=0$ for each arc $a \in A$,

\noindent for otherwise, we may decrease both $f(a)$ and $f(\bar{a})$ by their minimum.  

The above definition of $(S,T)$-flow is given in the arc-vertex form. An equivalent definition is formulated in the following path packing form: 
Let ${\cal P}$ be the set of all simple paths in $G$ from $S$ to $T$. An $(S,T)$-flow is an assignment $g:\, {\cal P}\rightarrow \mathbb R_+$ 
such that $\sum_{e \in Q \in {\cal P}}\, g(Q) \le  c(e)$ for all $e \in E$, whose value is $\sum_{Q \in {\cal P}}\, g(Q)$. We shall use both forms
in the subsequent sections.

As is well known, the following general max-flow min-cut theorem holds.

$(3.5)$ The maximum value of an $(S,T)$-flow is equal to the minimum capacity of an $(S,T)$-cut in the network $(G,c)$.

Given four distinguished vertices, two sources $s_1, s_2$ and two sinks $t_1,t_2$, of $G$, an $(s_1,t_1; s_2,t_2)$-{\em biflow} (or simply a {\em biflow})
consists of an $s_1$-$t_1$ flow $f_1$ and an $s_2$-$t_2$ flow $f_2$ in $G$, denoted by $(f_1,f_2)$, such that

$(3.6)$ $f_1(a)+f_1(\bar{a})+f_2(a)+f_2(\bar{a}) \le c(a)$ for each arc $a \in A$, 

\noindent whose {\em value} is defined to be ${\rm val}(f_1)+{\rm val}(f_2)$. Again, biflow in this arc-vertex form is equivalent to that in 
its path packing form. By definition, an $(s_1,t_1; s_2,t_2)$-{\em bicut} (or simply a {\em bicut}) is a subset $K$ of $E$ such that $G\de K$ contains
neither a path from $s_1$ to $t_1$ nor a path from $s_2$ to $t_2$. Hu \cite{Hu} proved that the minimum bicut can be found in polynomial time by
applying a maximum flow algorithm twice.

\begin{lemma}[Hu \cite{Hu}]\label{tchu}
The minimum capacity of an $(s_1,t_1; s_2,t_2)$-bicut is equal to the minimum capacity of an $(\{s_1,s_2\}, \{t_1,t_2\})$- or 
$(\{s_1,t_2\}, \{t_1,s_2\})$-cut. 
\end{lemma}

The following lemma will be used repeatedly in our design of an efficient algorithm for finding maximum integral biflows 
in Seymour graphs.  

\begin{lemma}\label{triflow}
Let $G=(V,E)$ be a graph with three distinguished vertices $x,y,z$, let $c:\, E \rightarrow \mathbb Z_+$ be a capacity function, and let $\tau_x$ (resp.
$\tau_y$) be the minimum capacity of an $(x, \{z,y\})$-cut (resp. $(\{x,z\},y)$-cut). Then the minimum capacity $\tau$ of an $x$-$y$ cut is equal to 
$\min\{\tau_x, \tau_y\}$. Furthermore, there exist
\begin{itemize}
\vspace{-2mm}
\item[(i)] an integral $(x, \{z,y\})$-flow $f$, an integral $x$-$y$ flow $f_1$, and an integral $x$-$z$ flow $f_2$ in the network $(D,c)$, such that
${\rm val}(f)={\rm div}_f(x)=\tau_x$, ${\rm val}(f_1)=-{\rm div}_f(y)=\tau$, ${\rm val}(f_2)=-{\rm div}_f(z)=\tau_x-\tau$, and
$f(a)=f_1(a)+f_2(a)$ for each arc $a$ in $D$;  and
\vspace{-2mm}
\item[(ii)] an integral $(\{x,z\},y)$-flow $g$, an integral $x$-$y$ flow $g_1$, and an integral $z$-$y$ flow $g_2$ in the network $(D,c)$, such that
${\rm val}(g)=-{\rm div}_g(y)=\tau_y$, ${\rm val}(g_1)={\rm div}_g(x)=\tau$, ${\rm val}(g_2)={\rm div}_g(z)=\tau_y-\tau$, and
$g(a)=g_1(a)+g_2(a)$ for each arc $a$ in $D$,
\vspace{-1mm}
\end{itemize}
\noindent all of which can be found in $O(nm)$ time, where $n=|V|$ and $m=|E|$. 

\end{lemma}

{\bf Proof.} Trivially, $\tau=\min\{\tau_x, \tau_y\}$. It suffices to prove statement (i), as statement (ii) is symmetric to (i).
In our proof, $G^*$ is the graph obtained from $G$ by adding a sink $t$ and two edges $yt$ and $zt$. For each subgraph $H$ of $G$ and $U \subseteq V(H)$, 
we use $\delta_H(U)$ to denote the set of all edges of $H$ between $U$ and $V(H)\de U$. Let us show that

(1) the minimum capacity of a $(z, \{x,y\})$-cut in $G$ is at least $\tau_x-\tau$.

For this purpose, let $\delta_G(J)$ be a minimum $(z, \{x,y\})$-cut in $G$ with $z\in J$, and let $H:=G\de J$. For any minimum $(x,y)$-cut $\delta_G(K)$ 
in $G$, clearly $\delta_G(K) \cap E(H)$ is an $x$-$y$ cut in $H$. So the minimum capacity of an $x$-$y$ cut $\delta_H(L)$ (with $x\in L$) 
in $H$ is at most $\tau$.  Note that $\delta_G(L)$ is an $(x, \{z,y\})$-cut in $G$ with capacity at most $c(\delta_H(L))+c(\delta_G(J)) 
\le \tau+c(\delta_G(J))$. So $\tau_x \le \tau+c(\delta_G(J))$, which implies $c(\delta_G(J))\ge \tau_x-\tau$, as desired. 

Let us first find the minimum capacity $\tau$ of an $x$-$y$ cut in $(G,c)$ using Orlin's algorithm \cite{O} along with King, Rao and Tarjan's algorithm \cite{KRT},
which takes $O(nm)$ time. To find the minimum capacity $\tau_x$ of an $(x, \{z,y\})$-cut, define $c(yt)=c(zt):=\infty$. Then a minimum $x$-$t$ cut in $(G^*,c)$ is 
a minimum $(x, \{z,y\})$-cut in $(G, c)$. So $T_x$ can also be obtained in $O(nm)$ time.

To find the flow $f$ as specified in (i), redefine $c(yt):=\tau$ and $c(zt):=\tau_x-\tau$. From (1) we deduce that the minimum capacity of an $x$-$t$ cut 
in $G^*$ is exactly $\tau_x$. By the max-flow min-cut theorem, $G^*$ contains a maximum integral $x$-$t$ flow $h$ with value $\tau_x$; the restriction of 
$h$ to $G$ is an integral $(x, \{z,y\})$-flow $f$ in $G$ with ${\rm val}(f)={\rm div}_f(x)=\tau_x$, ${\rm div}_f(y)=-\tau$, and ${\rm div}_f(z)=\tau-\tau_x$. 
The flow $h$ and hence $f$ can be obtained in $O(nm)$ time using the aforementioned algorithms. 

To find the flow $f_1$ as described in (i), let $D'$ be obtained from $D$ by adding two arcs $(y,t)$ and $(z,t)$, and let $c'$ be the capacity function defined 
on $D'$ by: $c'(a):=f(a)$ for each $a$ in $D$, $c'(y,t):=\tau$, and $c'(z,t):=0$. Then the minimum capacity of an $x$-$t$ cut in $(D', c')$ is exactly $\tau$. So 
we can find a maximum $x$-$t$ flow $l$ in $(D', c')$ with value $\tau$ in $O(nm)$ time; the restriction of $l$ to $(D,c)$ is an integral $x$-$y$-flow $f_1$ in 
$(D,c)$ with ${\rm val}(f_1)=-{\rm div}_f(y)=\tau$. Define $f_2(a):=f(a)-f_1(a)$ for each arc $a$ in $D$. From (3.1)-(3.4) and (3.6), we see that $f_2$ is an 
integral $x$-$z$ flow in the network $(D,c)$, such that ${\rm val}(f_2)=-{\rm div}_f(z)=\tau_x-\tau$. \qed

\section{Algorithm I}

In this section the Seymour graph $G=(V,E)$ we consider is as described in Theorem \ref{Thm3}(i); that is,

(4.1)  $G$ arises from a plane graph $H$ in $\mathcal{G}_1$ by 3-summing 3-connected graphs to its facial triangles, in which $s_1,s_2,t_2,t_1$ occur on 
the facial cycle $C$ of $H$ in the clockwise order.

Our objective is to design a combinatorial polynomial-time algorithm for finding a maximum integral $(s_1,t_1; s_2,t_2)$-biflow in the network $(G,c)$, where 
$c:\, E \rightarrow \mathbb Z_+$ is a capacity function.  

Although $G$ itself is not necessarily a planar graph, it inherits from its base graph $H$ some planarity properties, as shown below. 

\begin{lemma}\label{cross}
If $P_1$ is an $s_1$-$t_2$ path and $P_2$ is an $s_2$-$t_1$ path in $G$, then $P_1$ and $P_2$ have at least one internal vertex in common. 
\end{lemma}

\vspace{-1mm}
{\bf Proof.}  Let $a_1,a_2,\ldots, a_k$ be all the vertices in $V(H) \cap V(P_1)$ that occur on $P_1$ in order as we traverse $P_1$ from $s_1$, where $a_1=s_1$
and $a_k=t_2$. For each $1\le i \le k-1$,  the subpath $P[a_i, a_{i+1}]$ is either an edge $a_ia_{i+1}$ of $H$ or a path contained in a $3$-connected graph  
that is $3$-summed to a facial triangle $T_i$ of $H$, with $\{a_i,a_{i+1}\} \subseteq V(T_i)$ (so $a_ia_{i+1}\in  E(H)$).   
Thus $P_1'=a_1a_2\ldots a_k$ is a path connecting $s_1$ and $t_2$ in $H$. Similarly, we can define a path $P_2'$ connecting $s_2$ and $t_1$ in $H$, with
$V(P_2') \subseteq V(P_2)$. Since $s_1,s_2,t_2,t_1$ occur on the facial cycle $C$ of $H$ in the clockwise order, $P_1'$ and $P_2'$ (and hence $P_1$ and $P_2$)
must have at least one vertex in common. \qed

\vskip 4mm

\begin{lemma}\label{correspondence}
There is an integral $(s_1,t_1;s_2,t_2)$-biflow $(f_1,f_2)$ in the network $(G,c)$ with ${\rm val}(f_i)=k_i$ for $i=1,2$, if and only if there is an 
integral $(\{s_1,s_2\}, \{t_1,t_2\})$-flow $f$ in $(G,c)$ with ${\rm val}(f)=k_1+k_2$ and ${\rm div}_f(s_i)=-{\rm div}_f(t_i)=k_i$ for $i=1,2$. 
Furthermore, give such a biflow $(f_1,f_2)$, one can construct the desired flow $f$ in $O(m)$ time. Conversely, given such a flow $f$,  
one can construct the desired biflow $(f_1,f_2)$ in $O(m^2)$ time, where $m=|E(G)|$.
\end{lemma}

\vspace{-1mm}
{\bf Proof.} In our proof $D=(V,A)$ is the digraph obtained from $G$ by replacing each edge $uv \in E$ with a pair of opposite arcs $(u,v)$ and $(v,u)$, such 
that $c(u,v)=c(v,u)=c(uv)$. Recall that a flow in the undirected network $(G,c)$ is actually one in the directed network $(D,c)$. 

The ``ony if" part. Let $(f_1,f_2)$ be an integral $(s_1,t_1;s_2,t_2)$-biflow in the network $(G,c)$ with ${\rm val}(f_i)=k_i$ for $i=1,2$. 
Then ${\rm div}_{f_i}(s_i)=-{\rm div}_{f_i}(t_i)=k_i$ and ${\rm div}_{f_i}(v)=0$ for all vertex $v \ne s_i,t_i$ for $i=1,2$. Define $f(a)=f_1(a)+f_2(a)$ for 
each arc $a \in A$. Then ${\rm div}_f(v)= {\rm div}_{f_1}(v) +{\rm div}_{f_2}(v)$ for each vertex $v$. So ${\rm div}_f(v)=0$ for all $v\notin \{s_1,s_2,t_1,t_2\}$, 
${\rm div}_f(s_i)={\rm div}_{f_i}(s_i)=k_i$ and ${\rm div}_f(t_i)={\rm div}_{f_i}(t_i)=-k_i$ for $i=1,2$. From (3.1)-(3.5) we see that $f$ is an integral $(\{s_1,s_2\}, 
\{t_1,t_2\})$-flow in $(G,c)$ with value $k_1+k_2$ and ${\rm div}_f(s_i)=-{\rm div}_f(t_i)=k_i$ for $i=1,2$. Clearly, $f$ can be constructed in $O(m)$ time.

The ``if" part. Let $f$ be an integral $(\{s_1,s_2\}, \{t_1,t_2\})$-flow $f$ in the network $(G,c)$ with value $k_1+k_2$ and ${\rm div}_f(s_i)=
-{\rm div}_f(t_i)=k_i$ for $i=1,2$. Using breadth-first search $O(m)$ times, we can construct a collection ${\cal P}$ of $O(m)$ directed paths $Q$ in $D$ 
from $\{s_1,s_2\}$ to $\{t_1,t_2\}$ together with integral multiplicities $g(Q)>0$, such that  

(1) $\sum_{a \in Q \in {\cal P}}\, g(Q) \le f(a)$ for each $a \in A$;
 
(2) $\sum_{Q \in {\cal P}}\, g(Q)=k_1+k_2$; and

(3) $\sum_{v \in Q \in {\cal P}}\, g(Q)=k_i$ for $v=s_i,t_i$ and $i=1,2$.

To justify this, find a directed path $Q$ from $S$ to $T$ in $D$ such that $f(a)>0$ for each arc $a$ on $Q$ using breadth-first search (see, for instance,
\cite{CL}), which runs in $O(n+m)$ time.  Define $g(Q):=\delta=\min \{f(a): a \in P\}$. Replacing $f(a)$ by $f(a)-\delta$ for each arc $a$ on $P$, and 
repeat the process. This step requires $O((n+m)m)=O(m^2)$ time. 

Let ${\cal P}_i$ be the collection of directed paths $Q$ from $s_i$ to $t_i$ in ${\cal P}$, and let ${\cal P}_{2+i}$ be the collection 
of directed paths $Q$ from $s_i$ to $t_{3-i}$ in ${\cal P}$ for $i=1,2$. From (3) we deduce that

(4) $\sum_{Q \in {\cal P}_3}\, g(Q) = \sum_{Q \in {\cal P}_4}\, g(Q)$.

So ${\cal P}_3\ne \emptyset$ if and only if ${\cal P}_4\ne \emptyset$. Assume ${\cal P}_3\ne \emptyset \ne {\cal P}_4$. Let $P$ (resp. $Q$) be a path in ${\cal P}_3$ 
(resp. ${\cal P}_4$) and set $\delta:=\min\{g(P), g(Q)\}$.  By Lemma \ref{cross}, $P$ and $Q$ have at least one internal vertex in common. So there exist 
two arc-disjoint directed paths $P'$ from $s_1$ to $t_1$ and $Q'$ from $s_2$ to $t_2$, with all arcs contained in $P\cup Q$.  Define $g(R)=0$ if $R$ is a path not 
contained in ${\cal P}$. Set 

$\bullet$ $g(R):=g(R)-\delta$ if $R\in \{P, Q\}$ and $g(R):=g(R)+\delta$ if $R\in \{P',Q'\}$,

$\bullet$ ${\cal P}_1:= {\cal P}_1 \cup \{P'\}$, ${\cal P}_2:= {\cal P}_2 \cup \{Q'\}$, and ${\cal P}:= {\cal P} \cup \{P',Q'\}$.  

\noindent Remember to remove all paths $R$ with $g(R)=0$ in ${\cal P}$ and ${\cal P}_i$ for $1\le i \le 4$. Repeating this uncrossing process, we shall end up with
a collection ${\cal P}$ of directed paths from $\{s_1,s_2\}$ to $\{t_1,t_2\}$ with ${\cal P}_3=\emptyset={\cal P}_4$ by (4). This step requires $O(nm)$ time.
Define

$\bullet$ $f_i(a)=\sum_{a \in Q \in {\cal P}_i}\, g(Q)$ for each $a \in A$ and $i=1,2$.

\noindent Clearly, $f_i$ is an integral $s_i$-$t_i$ flow with ${\rm val}(f_i)={\rm div}_f(s_i)=k_i$ by (3) and can be constructed in $O(m^2)$ time
for $i=1,2$ . It follows from (2) that  ${\rm val}(f_1) +{\rm val}(f_2)=k_1+k_2$. By virtue of (3.2) and (1), we see that $(f_1,f_2)$ is a biflow in the 
network $(G,c)$ with ${\rm val}(f_i)=k_i$ for $i=1,2$. \qed

\vskip 4mm

The algorithm given in this section will also serve as a subroutine for that to be devised in the succeeding section. So we present the algorithm in a slightly
more general setting. Recall that, by Lemma \ref{tchu}, the minimum bicut can be found in polynomial time using a maximum flow algorithm.

\begin{theorem}\label{case1}
Let $\tau$ be the minimum capacity of an $(s_1,t_1;s_2,t_2)$-bicut in the network $(G,c)$, let $\tau_i$ be the minimum capacity of an $s_i$-$t_i$ cut for $i=1,2$, and
let $k_1$ and $k_2$ be two nonnegative integers satisfying $k_1\leq \tau_1$, $k_2\leq \tau_2$ and $k_1+k_2\leq \tau$. Then an integral $(s_1,t_1;s_2,t_2)$-biflow $(f_1,f_2)$ 
in $(G,c)$ with ${\rm val}(f_1)=k_1$ and ${\rm val} (f_2)=k_2$ exists and can be found in $O(m^2)$ time, where $m=|E(G)|$.
\end{theorem}

\vspace{-1mm}
{\bf Proof.} Once again, let $D=(V,A)$ be the digraph obtained from $G$ by replacing each edge $uv \in E$ with a pair of opposite arcs $(u,v)$ and $(v,u)$, such 
that $c(u,v)=c(v,u)=c(uv)$. Recall that a flow in the undirected network $(G,c)$ is actually one in the directed network $(D,c)$. 

Let us first find an integral $s_1$-$t_1$ flow $g$ in $(G,c)$ with ${\rm val}(g)={\rm div}_g(s_1)=k_1$ using Orlin's algorithm \cite{O} together with King, Rao and Tarjan's
algorithm \cite{KRT} for the maximum flow problem, which takes $O(nm)$ time. (The problem of finding $g$ can be reduced to a maximum flow problem: add a pendent arc 
$(t,t')$ with capacity $k_1$ to $D$ and view $t'$ as the sink.) Define a new capacity function $c':\, A \rightarrow \mathbb Z_+$ by

(1) $c'(a)=c(a)-g(a)+g(\bar{a})$ for each $a\in A$, where $\bar{a}$ is the reverse of $a$ in $D$. 

For a subset $U\subseteq V$, we have $c'(\delta^+(U))=\sum_{a\in \delta^+(U)} c'(a)= \sum_{a\in \delta^+(U)} [c(a)-g(a)+g(\bar{a})]
= c(\delta^+(U))-\sum_{a\in \delta^+(U)} [g(a)-g(\bar{a})]$. So 

(2) $c'(\delta^+(U)) = c(\delta^+(U))-\sum_{u\in U} {\rm div}_g(u)$ for any $U \subseteq V$. 
    
(3) $c'(\delta^+(U)) \ge k_2$	for any $U \subseteq V$ with $s_2\in U$ and $t_2 \notin U$.	

To justify this, we consider three cases, depending on whether $s_1$ or $t_1$ is contained in $U$.

$\bullet$ $s_1\in U$ and $t_1\notin U$. In this case, $\delta(U)$ is an $(s_1,t_1;s_2,t_2)$-bicut in the network $(G,c)$. So $c(\delta(U)) \ge \tau$.
Note that ${\rm div}_g(u)=0$ for any $u \in U\de \{s_1\}$. By (2), we obtain $c'(\delta^+(U)) = c(\delta^+(U))-\sum_{u\in U} {\rm div}_g(u)\ge
\tau- {\rm div}_g(s_1) =\tau-k_1\ge k_2$.  

$\bullet$ $s_1\notin U$ and $t_1\in U$. In this case, $\delta(U)$ is an $(\{s_1,t_2\}, \{t_1,s_2\})$-cut in the network $(G,c)$. By  
Lemma \ref{tchu}, we have $c(\delta(U)) \ge \tau$. Note that ${\rm div}_g(u)=0$ for any $u \in U\de \{t_1\}$. By (2), we obtain 
$c'(\delta^+(U)) = c(\delta^+(U))-\sum_{u\in U} {\rm div}_g(u)\ge \tau- {\rm div}_g(t_1) =\tau+ k_1\ge k_2$.  

$\bullet$ $\{s_1,t_1\} \subseteq U$ or $\{s_1, t_1\} \subseteq V\de U$. In this case, $\sum_{u\in U} {\rm div}_g(u)=0$. Since $\delta(U)$
is an $s_2$-$t_2$ cut in $(G,c)$, we have $c(\delta(U))\ge \tau_2 \ge k_2$. It follows from (2) that $c'(\delta^+(U)) = c(\delta^+(U))-
\sum_{u\in U} {\rm div}_g(u)\ge k_2$. 

Combining the above three cases, we see that (3) holds. 

By (3) and the max-flow min-cut theorem, we can then find an integral $s_2$-$t_2$ flow $h$ in the network $(G,c')$ with ${\rm val}(h)={\rm div}_h(s_2)=k_2$
in $O(nm)$ time using the aforementioned algorithms (see the technique for finding $g$ at the beginning of our proof). Define $f:\, A \rightarrow 
\mathbb Z_+$ by

(4) $f(a)=\max\{0,g(a)+h(a)-g(\bar{a})-h(\bar{a})\}$ for each $a\in A$. 

Observe that, for each $a\in A$, we have $h(a) \le c'(a)=c(a)-g(a)+g(\bar{a})$ by (1). So $g(a)+h(a)-g(\bar{a}) \le c(a)$ and hence $f(a) \le c(a)$. 
By definition, $f(a)f(\bar{a})=0$, which yields $f(a)+f(\bar{a}) \le c(a)$. For each $v \in V$, by definition ${\rm div}_f(v)=f(\delta^+(v))-f(\delta^-(v))=
\sum_{a\in\delta_{}^+(v)} (f(a) -f(\bar{a}))=\sum_{a\in\delta_{}^+(v)}(g(a)+h(a)-g(\bar{a})-h(\bar{a}))={\rm div}_g(v)+{\rm div}_h(v)$, where 
the third equality holds because $\max\{0, x-y\}- \max\{0, y-x\}=x-y$ for any real numbers $x$ and $y$. So ${\rm div}_f(s_1)={\rm div}_g(s_1)
=-{\rm div}_g(t_1)=-{\rm div}_f(t_1)$ and ${\rm div}_f(s_2)={\rm div}_h(s_2)=-{\rm div}_h(t_2)=-{\rm div}_f(t_2)$. Furthermore, 
${\rm div}_f(v)=0$ if $v \notin \{s_1,t_1,s_2,t_2\}$.  It follows from (3.1)-(3.3) that $f$ is an $(\{s_1,t_2\}, \{t_1,s_2\})$-flow with value
${\rm div}_f(s_1)+{\rm div}_f(s_2)={\rm div}_g(s_1)+{\rm div}_h(s_2)=k_1+k_2$. Therefore, Lemma \ref{correspondence} guarantees the existence 
of an integral $(s_1,t_1;s_2,t_2)$-biflow $(f_1,f_2)$ with ${\rm val}(f_1)=k_1$ and ${\rm val} (f_2)=k_2$ in the network $(G,c)$, which can 
be found in $O(m^2)$ time. \qed

\section{Algorithm II}

In this section the Seymour graph $G=(V,E)$ we consider is as described in Theorem \ref{Thm3}(iii); that is,

(5.1) $G$ arises from a graph $H'$ in $\mathcal{G}_3$ by 3-summing 3-connected graphs to triangles $s_1t_1u_1$, $s_1t_1v_1$, $s_2t_2u_2$, $s_2t_2v_2$, 
and to facial triangles of $H$. 

Our objective is to design a combinatorial polynomial-time algorithm for finding a maximum integral $(s_1,t_1; s_2,t_2)$-biflow in the network $(G,c)$, where 
$c:\, E \rightarrow \mathbb Z_+$ is a capacity function.

(5.2) The case when a Seymour graph  $G=(V,E)$ satisfies Theorem \ref{Thm3}(ii) is contained in (5.1) as a subcase.

To see this, recall that now $G$ arises from a graph in $\mathcal{G}_2$ by 3-summing 3-connected graphs to triangles $s_1t_1u_1$, $s_1t_1v_1$, and to facial triangles of $H$.
Symmetry allows us to assume that $u_1, v_1, e_2$ are on the facial cycle $C$ of $H$ (see Figure \ref{fig:graph in G2}).  Let $G'$ be obtained from $G$ by adding two new vertices
$s_2', t_2'$ and five edges $s_2't_2', s_2's_2, s_2't_2, t_2's_2, t_2't_2$, and view $e_1$ and $e_2'=s_2't_2'$ (rather than $e_2$) as two specified edges 
in $G'$. Let $c'$ be the capacity function obtained from $c$ by replacing $c(e_2)$ with $c'(e_2)=0$ and defining $c'(e_2')=c(e_2)$, $c'(s_2's_2)=c'(t_2't_2)
=\infty$, and $c'(s_2't_2)=c'(t_2's_2)=0$. Then there is one-to-one correspondence between integral $(s_1,t_1; s_2,t_2)$-biflows in the network $(G,c)$ 
and integral $(s_1,t_1; s_2',t_2')$-biflows in the network $(G',c')$. Since $G'$ is as described in Theorem \ref{Thm3}(iii), we get (5.2).

\begin{theorem}\label{case3}
A maximum integral $(s_1,t_1;s_2,t_2)$-biflow in the network $(G,c)$ can be found in $O(m^2)$ time, where $m=|E(G)|$.
\end{theorem}

\vspace{-1mm}
{\bf Proof.} By the definition of $\mathcal{G}_3$, the graph $H'$ in (5.1) is obtained from the complete graph $H$ on two vertices $u_1, v_1$ (now define $u_2=u_1$ 
and $v_2=v_1$) or from a $2$-connected plane graph $H$ with vertices $v_1,u_1, u_2, v_2$ on its outer facial cycle in the clockwise order by adding four vertices 
$s_1, t_1, s_2, t_2$ and ten edges $s_1t_1, s_1u_1, s_1v_1, t_1u_1, t_1v_1, s_2t_2, s_2u_2, s_2v_2, t_2u_2, t_2v_2$, where $u_i\ne v_i$ for $i=1,2$ and 
$|\{u_1,v_1, u_2,v_2\}|\ge 3$ (see Figure \ref{fig:graph in G3}; possibly $\{u_1,v_1\} \cap \{u_2,v_2\} \ne \emptyset$).   

Observe that $G$ is the union of the following five subgraphs $G_0, G_1, \ldots, G_4$: 

$\bullet$ $G_0$ is $H$ if $H$ is a complete graph on two vertices $u_1, v_1$  and arises from the plane graph $H$ by 3-summing 3-connected graphs to 
its facial triangles otherwise;

$\bullet$ $G_1$ arises from the triangle $T_1=s_1t_1u_1$ by possibly 3-summing a 3-connected graph to it; 

$\bullet$ $G_2$ arises from the triangle $T_2=s_1t_1v_1$ by possibly 3-summing a 3-connected graph to it; 

$\bullet$ $G_3$ arises from the triangle $T_3=s_2t_2u_2$ by possibly 3-summing a 3-connected graph to it; 

$\bullet$ $G_4$ arises from the triangle $T_4=s_2t_2v_2$ by possibly 3-summing a 3-connected graph to it. 

To find a maximum integral $(s_1,t_1;s_2,t_2)$-biflow in the network $(G,c)$, we shall turn to finding an appropriate flow in each network 
$(G_i, c_i)$ for $0\le i \le 4$, where the capacity function $c_i$ is the restriction of $c$ to $G_i$, and eventually combine them into a desired biflow.    
For convenience, we assume that

(1) $c(e_1)=c(e_2)=0$.

Otherwise, let $c':\, E \rightarrow \mathbb Z_+$ be obtained from the capacity function $c$ by replacing each $c(e_i)$ with $0$ for $i=1,2$. Find a maximum 
integral $(s_1,t_1; s_2,t_2)$-biflow $(f_1',f_2')$ in the network $(G,c')$. Let $f_i$ be obtained from $f_i'$ be setting $f_i(e_i)=c(e_i)$ for $i=1,2$. Clearly, 
$(f_1,f_2)$ is a maximum integral $(s_1,t_1; s_2,t_2)$-biflow in the network $(G,c)$. So a general capacity function $c$ can be reduced to one satisfying (1). 

We use the same notation as introduced in Lemma \ref{triflow}. 

(2) Let $\tau_{s_1}^i$ (resp. $\tau_{t_1}^i$) be the capacity of a minimum $(s_1, V(T_i)\de s_1)$-cut $K_{s_1}^i$ (resp. $(V(T_i)\de t_1, t_1)$-cut $K_{t_1}^i$) 
in the network $(G_i,c_i)$, and let $\tau^i$ be the capacity of a minimum $s_1$-$t_1$ cut in the network $(G_i,c_i)$ for $i=1,2$. 
The equality $\tau^i=\min\{\tau_{s_1}^i,\tau_{t_1}^i\}$ holds.

(3) Let $\tau_{s_2}^i$ (resp. $\tau_{t_2}^i$) be the capacity of a minimum $(s_2, V(T_i)\de s_2)$-cut $K_{s_2}^i$ (resp. $(V(T_i)\de t_2, t_2)$-cut $K_{t_2}^i$) 
in the network $(G_i,c_i)$, and let $\tau^i$ be the capacity of a minimum $s_2$-$t_2$ cut in the network $(G_i,c_i)$ for $i=3,4$.
The equality $\tau^i=\min\{\tau_{s_2}^i,\tau_{t_2}^i\}$ holds.

We say that $G_1$ and $G_2$ are {\em consistent} if there is a permutation $\{a,b\}$ of $\{s_1,t_1\}$ such that $\tau_{a}^i \ge \tau_{b}^i$ for $i=1,2$
and {\em inconsistent} otherwise. Similarly, we can define consistency and inconsistency for $G_3$ and $G_4$. We proceed by considering two cases.

{\bf Case 1.} $G_j$ and $G_{j+1}$ are consistent for $j=1$ or 3. 

In this case, by symmetry we may assume that $j=1$ and $\tau_{s_1}^i \ge \tau_{t_1}^i$ for $i=1,2$. Thus Lemma \ref{triflow} guarantees the existence of 

(4) an integral $(s_1, \{u_1,t_1\})$-flow $f_1$ in $G_1$ with ${\rm div}_{f_1}(s_1)=\tau_{s_1}^1$, ${\rm div}_{f_1}(u_1)=\tau^1-\tau_{s_1}^1$, and 
${\rm div}_{f_1}(t_1)=-\tau^1$ (see (2)) and  
 
(5) an integral $(s_1, \{v_1,t_1\})$-flow $f_2$ in $G_2$ with ${\rm div}_{f_2}(s_1)=\tau_{s_1}^2$, ${\rm div}_{f_2}(v_1)=\tau^2-\tau_{s_1}^2$, and 
${\rm div}_{f_2}(t_1)=-\tau^2$ (see (2)).   
 
Let us decompose $f_1$ into an integral $s_1$-$t_1$ flow $g_1$ and an integral $s_1$-$u_1$ flow $h_1$ in $G_1$ with ${\rm val}(g_1)=\tau^1$ 
and ${\rm val}(h_1)=\tau_{s_1}^1-\tau^1$, and decompose $f_2$ into an integral $s_1$-$t_1$ flow $g_2$ and an integral $s_1$-$v_1$ flow $h_2$ in 
$G_2$ with ${\rm val}(g_2)=\tau^2$ and ${\rm val}(h_2)=\tau_{s_1}^2-\tau^2$; such decompositions are available by (4) and (5).  
Define $z_1=g_1+g_2$. Then 

(6) $z_1$ is an $s_1$-$t_1$ flow in $G_1 \cup G_2$ and hence in $G$ with value $\tau^1+\tau^2$. 

Set $\theta:=\min \{{\rm val}(h_1), {\rm val}(h_2) \} =\min \{\tau_{s_1}^1-\tau^1, \tau_{s_1}^2-\tau^2\}$.
Let $h_i^*$ be the flow obtained from $h_i$ by reversing the direction of each arc flow for $i=1,2$. Then we can concatenate $h_1^*$ and $h_2$ to form
an integral $u_1$-$v_1$ flow $h_{12}$ with value $\theta$, and concatenate $h_2^*$ and $h_1$ to form an integral $v_1$-$u_1$ flow $h_{21}$ with value $\theta$.  
 
Let $G'$ be the union of $G_0 \cup G_3 \cup G_4$ and the edge $u_1v_1$, and let $c'$ be obtained from the restriction of $c$ to $G_0 \cup G_3 \cup G_4$ by
defining $c'(u_1v_1)=c(u_1v_1)+\theta$, where $c(u_1v_1)=0$ if $u_1v_1 \notin E(G)$. Using a maximum flow algorithm, find a 
maximum integral $s_2$-$t_2$ flow $f'$ and a minimum $s_2$-$t_2$ cut $K'$ in $(G', c')$ (actually in $(D', c')$, where $D'$ is obtained from $G'$ be 
replacing each arc with a pair of opposite arcs). By (3.4), we have 

(7) $f'(u_1,v_1)f'(v_1,u_1)=0$.  

\noindent The max-flow min-cut theorem also yields

(8) $f'(u_1,v_1)=c'(u_1,v_1)$ if $(u_1,v_1) \in K'$ and $f'(v_1,u_1)=c'(v_1,u_1)$ if $(v_1,u_1) \in K'$.

By replacing certain units $(\le \theta)$ of the flow on arc $(u_1, v_1)$ or $(v_1,u_1)$ (see (7)) in $f'$ with a subflow of
$h_{12}$ or $h_{21}$, we can get an $s_2$-$t_2$ flow $z_2$ in $G$ with the same value as $f'$. 

(9) $(z_1,z_2)$ is a maximum integral $(s_1,t_1;s_2,t_2)$-biflow in the network $(G,c)$.

To justify this, we consider two subcases.

$\bullet$ $u_1v_1 \notin K'$. In this subcase, $K'$ is also an $s_2$-$t_2$ cut in $G$, for otherwise, $G\de K'$ would contain an $s_2$-$t_2$ path $P$. 
Since $K'$ is an $s_2$-$t_2$ cut in $G'$, path $P$ contains at least one edge of $G_1\cup G_2$. So it contains both $u_1$ and $v_1$. Let $P'$ be obtained 
from $P$ by replacing $P[u_1,v_1]$ with $u_1v_1$. Then $P'$ is an $s_2$-$t_2$ path in $G'$ that contains no edge in $K'$, a contradiction. Let 
$K:=K_{t_1}^1 \cup K_{t_1}^2 \cup K'$ (see (2)). Then $K$ is an $(s_1,t_1;s_2,t_2)$-bicut in $(G,c)$, with capacity $c(K)= c(K_{t_1}^1)+c(K_{t_1}^2)+c(K')
=\tau^1+\tau^2+c'(K') =\tau^1+\tau^2+{\rm val}(z_2)$ (see (1)). It follows from (6) that $c(K)={\rm val}(z_1) +{\rm val}(z_2)$. Therefore, $K$ is a 
minimum $(s_1,t_1;s_2,t_2)$-bicut in the network $(G,c)$ and $(z_1,z_2)$ is a maximum integral $(s_1,t_1;s_2,t_2)$-biflow. 

$\bullet$ $u_1v_1 \in K'$. In this subcase, $c'(K')=c(K')+\theta$ by the definition of $c'(u_1v_1)$. Define $K''=K'$ if $u_1v_1 \in E(G)$ and $K'\de \{u_1v_1\}$ 
otherwise. Then $c(K'')=c(K')$. Symmetry allows us to assume that $\theta={\rm val}(h_1) =\tau_{s_1}^1-\tau^1$. Let $K=K_{s_1}^1 \cup K_{t_1}^2 \cup K''$ 
(see (2)). Clearly, $K$ is an $(s_1,t_1;s_2,t_2)$-bicut in $(G,c)$, with capacity $c(K)= c(K_{s_1}^1) +c(K_{t_1}^2)+ c(K'')=\tau^1_{s_1}+\tau^2+c(K')$ (see (1)). 
So $c(K)=(\tau^1+\theta)+\tau^2+c(K')=\tau^1+\tau^2+c'(K')=\tau^1+\tau^2+{\rm val}(z_2)$. It follows from (6) that $c(K)={\rm val}(z_1) +{\rm val}(z_2)$. Therefore, 
$K$ is a minimum $(s_1,t_1;s_2,t_2)$-bicut in the network $(G,c)$ and $(z_1,z_2)$ is a maximum integral $(s_1,t_1;s_2,t_2)$-biflow. So (9) is established.

{\bf Case 2.} $G_j$ and $G_{j+1}$ are inconsistent for $j=1$ and $3$. 

In this case, by symmetry we may assume that $\tau_{s_1}^1 \ge \tau_{t_1}^1$, $\tau_{t_1}^2 \ge \tau_{s_1}^2$, $\tau_{s_2}^3 \ge \tau_{t_2}^3$, and 
$\tau_{t_2}^4 \ge \tau_{s_2}^4$. Thus Lemma \ref{triflow} guarantees the existence of 

(10) an integral $(s_1, \{u_1,t_1\})$-flow $f_1$ in $G_1$ with ${\rm div}_{f_1}(s_1)=\tau_{s_1}^1$, ${\rm div}_{f_1}(u_1)=\tau^1-\tau_{s_1}^1$, and 
${\rm div}_{f_1}(t_1)=-\tau^1$ (see (2));  
 
(11) an integral $(\{s_1,v_1\},t_1)$-flow $f_2$ in $G_2$ with ${\rm div}_{f_2}(s_1)=\tau^2$, ${\rm div}_{f_2}(v_1)=\tau_{t_1}^2-\tau^2$, and 
${\rm div}_{f_2}(t_1)=-\tau^2_{t_1}$ (see (2));   

(12) an integral $(s_2, \{u_2,t_2\})$-flow $f_3$ in $G_3$ with ${\rm div}_{f_3}(s_2)=\tau_{s_2}^3$, ${\rm div}_{f_3}(u_2)=\tau^3-\tau_{s_2}^3$, and 
${\rm div}_{f_3}(t_2)=-\tau^3$ (see (3));  and
 
(13) an integral $(\{s_2,v_2\},t_2)$-flow $f_4$ in $G_4$ with ${\rm div}_{f_4}(s_2)=\tau^4$, ${\rm div}_{f_4}(v_2)=\tau_{t_2}^4-\tau^4$, and 
${\rm div}_{f_4}(t_2)=-\tau^4_{t_2}$ (see (3)).   

Let us decompose $f_1$ into an $s_1$-$t_1$ flow $g_1$ and an $s_1$-$u_1$ flow $h_1$ in $G_1$ with ${\rm val}(g_1)=\tau^1$ 
and ${\rm val}(h_1)=\tau_{s_1}^1-\tau^1$, decompose $f_2$ into an $s_1$-$t_1$ flow $g_2$ and an $v_1$-$t_1$ flow $h_2$ in 
$G_2$ with ${\rm val}(g_2)=\tau^2$ and ${\rm val}(h_2)=\tau_{t_1}^2-\tau^2$, decompose $f_3$ into an $s_2$-$t_2$ flow $g_3$ and 
an $s_2$-$u_2$ flow $h_3$ in $G_3$ with ${\rm val}(g_3)=\tau^3$ and ${\rm val}(h_3)=\tau_{s_2}^3-\tau^3$, and decompose $f_4$ into 
an $s_2$-$t_2$ flow $g_4$ and an $v_2$-$t_2$ flow $h_4$ in $G_4$ with ${\rm val}(g_4)=\tau^4$ and ${\rm val}(h_4)=\tau_{t_2}^4-\tau^4$,
where all $g_i$ and $h_i$ are integral; such decompositions are available by (10)-(13).  Set $\theta_1:=\min \{{\rm val}(h_1), {\rm val}(h_2) \} 
=\min \{\tau_{s_1}^1-\tau^1, \tau_{t_1}^2-\tau^2\}$ and $\theta_2:=\min \{{\rm val}(h_3), {\rm val}(h_4) \} =\min \{\tau_{s_2}^3-\tau^3, \tau_{t_2}^4-\tau^4\}$.

Let $\lambda$ be the minimum capacity of a $(u_1,v_1;u_2,v_2)$-bicut in $(G_0,c_0)$, and let $\lambda_i$ be the minimum capacity of a $u_i$-$v_i$ cut in $(G_0,c_0)$ for
$i=1,2$. Set $k_1:=\min\{\theta_1,\lambda_1\}$ and $k_2:=\min\{\theta_2, \lambda_2, \lambda-k_1\}$. By Theorem \ref{case1}, there exists a $(u_1,v_1;u_2,v_2)$-biflow 
$(r_1,r_2)$ in $(G_0,c_0)$ with ${\rm val}(r_1)=k_1$ and ${\rm val}(r_2)=k_2$. Thus we can concatenate $h_1$, $r_1$ and $h_2$ to form an integral $s_1$-$t_1$ flow $l_1$
with value $k_1$, and concatenate $h_3$, $r_2$ and $h_4$ to form an integral $s_2$-$t_2$ flow $l_2$ with value $k_2$.  Define $z_1=g_1+g_2+l_1$ and $z_2=g_3+g_4+l_2$.

(14) $(z_1,z_2)$ is a maximum integral $(s_1,t_1;s_2,t_2)$-biflow in the network $(G,c)$.

To justify this, we consider five subcases.

$\bullet$ $k_1=\theta_1$ and $k_2=\theta_2$. In this subcase, symmetry allows us to assume that $\theta_1={\rm val}(h_1) =\tau_{s_1}^1-\tau^1$ and
$\theta_2={\rm val}(h_3) =\tau_{s_2}^3-\tau^3$.  Let $K:=K_{s_1}^1 \cup K_{s_1}^2 \cup K_{s_2}^3 \cup K_{s_2}^4$ (see (2) and (3)). Then $K$ is an 
$(s_1,t_1;s_2,t_2)$-bicut in $(G,c)$, with capacity $c(K)= c(K_{s_1}^1)+ c(K_{s_1}^2)+ c(K_{s_2}^3)+ c(K_{s_2}^4) =\tau^1_{s_1}+\tau^2+
\tau^3_{s_2}+\tau^4=\tau^1+\tau^2+\tau^3+\tau^4+\theta_1+\theta_2=(\tau^1+\tau^2+k_1)+(\tau^3+\tau^4+k_2)$. So $c(K)={\rm val}(z_1) +{\rm val}(z_2)$. 
Therefore, $K$ is a minimum $(s_1,t_1;s_2,t_2)$-bicut in the network $(G,c)$ and $(z_1,z_2)$ is a maximum integral $(s_1,t_1;s_2,t_2)$-biflow. 

$\bullet$ $k_1=\theta_1$ and $k_2=\lambda_2$. In this subcase, symmetry allows us to assume that $\theta_1={\rm val}(h_1) =\tau_{s_1}^1-\tau^1$.
Let $C_2$ be a minimum $u_2$-$v_2$ cut in the network $(G_0,c_0)$, and let $K:=K_{s_1}^1 \cup K_{s_1}^2 \cup K_{t_2}^3 \cup K_{s_2}^4 \cup C_2$ (see (2) and (3)).
Then $K$ is an $(s_1,t_1;s_2,t_2)$-bicut in $(G,c)$, with capacity $c(K)= c(K_{s_1}^1)+ c(K_{s_1}^2)+ c(K_{t_2}^3)+ c(K_{s_2}^4) +c(C_2)=\tau^1_{s_1}+\tau^2+
\tau^3+\tau^4+ \lambda_2=\tau^1+\tau^2+\tau^3+\tau^4+\theta_1+\lambda_2=(\tau^1+\tau^2+k_1)+(\tau^3+\tau^4+k_2)$. So $c(K)={\rm val}(z_1) +{\rm val}(z_2)$. 
Therefore, $K$ is a minimum $(s_1,t_1;s_2,t_2)$-bicut in the network $(G,c)$ and $(z_1,z_2)$ is a maximum integral $(s_1,t_1;s_2,t_2)$-biflow. 

$\bullet$ $k_1=\lambda_1$ and $k_2=\theta_2$. This subcase is symmetric to the preceding one.

$\bullet$ $k_2=\lambda-k_1$.  In this subcase, let $C_0$ be a minimum $(u_1,v_1;u_2,v_2)$-bicut in the network $(G_0, c_0)$, and let 
$K:=K_{t_1}^1 \cup K_{s_1}^2 \cup K_{t_2}^3 \cup K_{s_2}^4 \cup C_0$. Then $K$ is an $(s_1,t_1;s_2,t_2)$-bicut in $(G,c)$, with capacity $c(K)= 
c(K_{t_1}^1)+ c(K_{s_1}^2)+ c(K_{t_2}^3)+ c(K_{s_2}^4) +c(C_0)  =\tau^1+\tau^2+ \tau^3+\tau^4+ \lambda = (\tau^1+\tau^2+k_1)+(\tau^3+\tau^4+k_2)$. 
So $c(K)={\rm val}(z_1) +{\rm val}(z_2)$. Therefore, $K$ is a minimum $(s_1,t_1;s_2,t_2)$-bicut in the network $(G,c)$ and $(z_1,z_2)$ is a maximum 
integral $(s_1,t_1;s_2,t_2)$-biflow. 

$\bullet$ $k_1=\lambda_1$ and $k_2=\lambda_2$. In this subcase, by the definition of $k_2$, we have $\lambda_2 \le \lambda-k_1 = \lambda- \lambda_1$.
So $\lambda \ge \lambda_1+\lambda_2$. Since the union of any $u_1$-$v_1$ cut and any $u_2$-$v_2$ cut separates these two pairs simultaneously, 
$\lambda \le \lambda_1+\lambda_2$. Hence $\lambda=\lambda_1+\lambda_2$, which implies $k_2=\lambda-k_1$. So the present case reduces to the 
preceding one. This proves (14).

In both cases, we need to apply Lemma \ref{triflow} once to each $G_i$ for $1 \le i \le 4$. In Case 1, the problem of constructing $h_{12}$ and $h_{21}$
are essentially maximum flow problems (see the technique for finding $g$ in the second paragraph of the proof of Theorem \ref{case1}). To find $f'$
(see the paragraph above (7)) and $z_2$ (see the paragraph below (8)), we need to apply a maximum flow algorithm respectively. In Case 2, to find 
$\lambda$, $\lambda_i$ for $i=1,2$, $l_1$, and $l_2$ (see the paragraph above (14)), we need to apply a maximum flow algorithm five times. To find
the biflow $(r_1,r_2)$, we need to apply Theorem \ref{case1} once. Thus our algorithm runs in $9\times O(nm)+O(m^2)=O(m^2)$ time. It is worthwhile
pointing out that cuts are only used to prove the optimality of the biflow constructed, so we need not output them explicitly. \qed

\section{Algorithm III}

In this section the Seymour graph $G=(V,E)$ we consider is as described in Theorem \ref{Thm3}(iv); that is,

(6.1) For each $M$-bridge $B$ of $G$, where $M$ is the graph formed by $e_1$ and $e_2$ only, either $B$ is trivial, or $B$ has exactly three feet, or $B$ has 
a 1-separation $(B_1,B_2)$ with $V(B_1)\cap V(B_2)=\{x\}$ such that $x\notin V(M)$ and $\{s_i,t_i\} \subseteq V(B_i)$ for $i=1,2$.

Our objective is to design a combinatorial polynomial-time algorithm for finding a maximum integral $(s_1,t_1; s_2,t_2)$-biflow in the network $(G,c)$, where 
$c:\, E \rightarrow \mathbb Z_+$ is a capacity function.  

\begin{theorem}\label{case4}
A maximum integral $(s_1,t_1;s_2,t_2)$-biflow in the network $(G,c)$ can be found in $O(nm)$ time, where $n=|V(G)|$ and $m=|E(G)|$.
\end{theorem}

\vspace{-1mm}
{\bf Proof.} For convenience, we assume that

(1) $c(e_1)=c(e_2)=0$ (see (1) in the proof of Theorem \ref{case3}). 

Let $B_1, B_2, \ldots, B_p$ be all $M$-bridges with exactly three feet, and let $B_{p+1}, B_{p+2}, \ldots, B_q$ be all $M$-bridges with four feet.
For $1\le i \le q$, let $c_i$ stand for the restriction of the capacity function $c$ to $B_i$. For $p+1\le i \le q$, let $(B_i^1,B_i^2)$ be a 1-separation 
of $B_i$ with $V(B_i^1)\cap V(B_i^2)=\{x_i\}$, such that $x_i\notin V(M)$ and $\{s_j,t_j\} \subseteq V(B_i^j)$ for $j=1,2$; with a slight abuse 
of notation, we still use $c_i$ to denote the restriction of the capacity function $c_i$ to $B_i^j$.  

We employ the same notation as introduced in Lemma \ref{triflow}. In the following statements (2)-(4), we assume that $V(B_i) \cap V(M)= \{s_j,t_j, r_i\}$ 
for $1\le i \le p$, where $j=1$ or $2$. 

(2) For $1\le i \le p$, let $\tau_{s_j}^i$ (resp. $\tau_{t_j}^i$) be the capacity of a minimum $(s_j, \{r_i,t_j\})$-cut $K_{s_j}^i$ (resp. $(\{s_j, r_i\},t_j)$-cut 
$K_{t_j}^i$) in the network $(B_i,c_i)$, and let $\tau^i$ be the capacity of a minimum $s_j$-$t_j$ cut in the network $(B_i,c_i)$. 
The equality $\tau^i=\min\{\tau_{s_j}^i,\tau_{t_j}^i\}$ holds.

Lemma \ref{triflow} guarantees the existence of an integral $s_j$-$t_j$ flow $f_i$ in the network $(B_i, c_i)$, as described below. 

(3) For $1\le i \le p$, if $\tau_{s_j}^i\ge \tau_{t_j}^i$, then $f_i$ is an integral $(s_j, \{r_i,t_j\})$-flow in $(B_i, c_i)$ with ${\rm div}_{f_i}(s_j)=\tau_{s_j}^i$, 
${\rm div}_{f_i}(r_i)=\tau^i-\tau_{s_j}^i$, and ${\rm div}_{f_i}(t_j)=-\tau^i$. We decompose $f_i$ into an integral $s_j$-$t_j$ flow $g_i$ and an integral $s_j$-$r_i$ 
flow $h_i$ in $(B_i, c_i)$ with ${\rm val}(g_i)=\tau^i$ and ${\rm val}(h_i)=\tau_{s_j}^i-\tau^i$. Let $J_i$ denote $K_{t_j}^i$ and $J_i'$ denote $K_{s_j}^i$.

(4) For $1\le i \le p$, if $\tau_{t_j}^i\ge \tau_{s_j}^i$, then $f_i$ is an integral $(\{s_j, r_i\}, t_j\})$-flow in $(B_i, c_i)$ with ${\rm div}_{f_i}(s_j)=\tau^i$, 
${\rm div}_{f_i}(r_i)=\tau_{t_j}^i-\tau^i$, and ${\rm div}_{f_i}(t_j)=- \tau_{t_j}^i$. We decompose $f_i$ into an integral $s_j$-$t_j$ flow $g_i$ and an integral 
$r_i$-$t_j$ flow $h_i$ in $(B_i, c_i)$ with ${\rm val}(g_i)=\tau^i$ and ${\rm val}(h_i)=\tau_{t_j}^i-\tau^i$. Let $J_i$ denote $K_{s_j}^i$ and $J_i'$ 
denote $K_{t_j}^i$.

We also need to find some integral flows in $B_i^1$ and $B_i^2$ for $p+1\le i \le q$ using Lemma \ref{triflow}. 

(5) For $p+1\le i \le q$ and $j=1,2$, let $\tau_{s_j}^i$ (resp. $\tau_{t_j}^i$) be the capacity of a minimum $(s_j, \{x_i,t_j\})$-cut $K_{s_j}^i$ (resp. 
$(\{s_j, x_i\},t_j)$-cut $K_{t_j}^i$) in the network $(B_i^j, c_i)$, and let $\tau^i_j$ be the capacity of a minimum $s_j$-$t_j$ cut in the network 
$(B_i^j,c_i)$. The equality $\tau^i_j=\min\{\tau_{s_j}^i,\tau_{t_j}^i\}$ holds.

(6) For $p+1\le i \le q$ and $j=1,2$, if $\tau_{s_j}^i\ge \tau_{t_j}^i$, then $f_{ij}$ is an integral $(s_j, \{x_i,t_j\})$-flow in $(B_i^j, c_i)$ with ${\rm div}_{f_{ij}}(s_j)
=\tau_{s_j}^i$, ${\rm div}_{f_{ij}}(x_i)=\tau^i_j-\tau_{s_j}^i$, and ${\rm div}_{f_{ij}}(t_j)=-\tau^i_j$. We decompose $f_{ij}$ into an integral $s_j$-$t_j$ flow $g_{ij}$ and an 
integral $s_j$-$x_i$ flow $h_{ij}$ in $(B_i^j, c_i)$ with ${\rm val}(g_{ij})=\tau^i_j$ and ${\rm val}(h_{ij})=\tau_{s_j}^i-\tau^i_j$. When $j=1$, let
$L_i$ denote $K_{t_1}^i$ and $L_i'$ denote $K_{s_1}^i$. When $j=2$, let $R_i$ denote $K_{t_2}^i$ and $R_i'$ denote $K_{s_2}^i$. 

(7) For $p+1\le i \le q$ and $j=1,2$, if $\tau_{t_j}^i\ge \tau_{s_j}^i$, then $f_{ij}$ is an integral $(\{s_j, x_i\},t_j)$-flow in $(B_i^j, c_i)$ with 
${\rm div}_{f_{ij}}(s_j) =\tau_j^i$, ${\rm div}_{f_{ij}}(x_i)=\tau_{t_j}^i-\tau^i_j$, and ${\rm div}_{f_{ij}}(t_j)=-\tau^i_{t_j}$. We decompose $f_{ij}$ into 
an integral $s_j$-$t_j$ flow $g_{ij}$ and an integral $x_i$-$t_j$ flow $h_{ij}$ in $(B_i^j, c_i)$ with ${\rm val}(g_{ij})=\tau^i_j$ and ${\rm val}(h_{ij})=
\tau_{t_j}^i-\tau^i_j$. When $j=1$, let $L_i$ denote $K_{s_1}^i$ and $L_i'$ denote $K_{t_1}^i$. When $j=2$, let $R_i$ denote $K_{s_2}^i$ and $R_i'$ denote 
$K_{t_2}^i$. 

For $p+1\le i \le q$, we can concatenate $h_{i1}$ and $h_{i2}$ defined in (6) and (7) to form an $(\{s_1, t_1\}, \{s_2, t_2\})$-flow $h_i$ in $(G,c)$ with
${\rm val}(h_i) = \min\{{\rm val}(h_{i1}), {\rm val}(h_{i2}) \}$. For each $a\in \{s_1, t_1\}$ and $b\in \{s_2, t_2\}$, let ${\Pi}_{ab}$ consist
of all subscripts $i$ with $1\le i \le q$ such that $h_i$ is an $a$-$b$ flow with ${\rm val}(h_i)>0$ (including those defined in (3) and (4)). 
 
Let us construct a complete bipartite graph $H$ on $V(M)$ with bipartition $(\{s_1, t_1\}, \{s_2, t_2\})$. For each edge $ab$ of $H$, define its capacity
$w(ab):=c(ab)+ \sum_{i \in \Pi_{ab}}  {\rm val}(h_i)$. Let $E_1:=\{s_1s_2,t_1t_2\}$ and $E_2:=\{s_1t_2,s_2t_1\}$. Then both $E_1$ and $E_2$ are
$(s_1, t_1; s_2, t_2)$-bicuts in $H$.

(8) A maximum integral $(s_1, t_1; s_2, t_2)$-biflow $(l_1, l_2)$ in $(H, w)$ with value $\min \{w(E_1), w(E_2)\}$ can be found in constant time.

To justify this, we first find a maximum integral $s_1$-$t_1$ flow $l_1$ and a minimum $s_1$-$t_1$ cut $C_1$ in $(H, w)$. If $C_1=E_i$ for $i=1$ 
or $2$, then $C_1$ is a minimum $(s_1, t_1; s_2, t_2)$-bicut and $(l_1, \emptyset)$ is a maximum integral $(s_1, t_1; s_2, t_2)$-biflow, 
we are done. So we assume that $C_1\ne E_i$ for $i=1,2$. From the structure of $H$, we deduce that $C_1$ consists of two edges incident with the same 
vertex $v$, say $v=s_1$. Thus $l_1(s_1s_2) =w(s_1s_2)=l_1(s_2t_1)$ and  $l_1(s_1t_2) =w(s_1t_2)=l_1(t_1t_2)$. Let $l_2$ be the $s_2$-$t_2$ flow of 
value $\delta:=\min \{w(s_2t_1)- w(s_1s_2), w(t_1t_2)- w(s_1t_2)\}$ shipped along the path $s_2t_1t_2$. Then ${\rm val}(l_1)+ {\rm val}(l_2)= w(s_1s_2) 
+w(s_1t_2)+\delta=\min \{w(E_1), w(E_2)\}$. So (8) holds.

For $i=1,2$, by replacing each arc flow $l_i(ab)$ with an integral $a$-$b$ flow contained in $\sum_{i\in \Pi_{ab}} h_i$ with the 
same value, we can get an $s_i$-$t_i$ flow $l_i^*$ in the network $(G,c)$ with ${\rm val}(l_i^*)={\rm val}(l_i)$. Besides, let
$\Lambda_i$ consist of all subscripts $j$ with $1 \le j \le p$ such that $g_j$ defined in (3) and (4) is an $s_i$-$t_i$ flow.  Define
$z_k:= \sum_{i \in \Lambda_k} g_i + \sum_{i=p+1}^q g_{ik}+l_k^*$ for $k=1,2$. Let us show that

(9) $(z_1,z_2)$ is a maximum integral $(s_1,t_1;s_2,t_2)$-biflow in the network $(G,c)$.

To this end, we assume, without loss of generality, that $\min \{w(E_1), w(E_2)\}=w(E_1)$. Let us divide the subscripts $1,2, \ldots, q$
into five sets:

$\bullet$ $\Gamma_1$ consists of all subscripts $i$ with $1\le i \le p$ in ${\Pi}_{s_1s_2} \cup  {\Pi}_{t_1t_2}$; 

$\bullet$ $\Gamma_2$ consists of all subscripts $i$ with $1\le i \le p$ outside $\Gamma_1$; 

$\bullet$ $\Gamma_3$ consists of all subscripts $i$ with $p+1\le i \le p$ in ${\Pi}_{s_1s_2} \cup  {\Pi}_{t_1t_2}$ with  ${\rm val}(h_{i1}) \le {\rm val}(h_{i2})$;

$\bullet$ $\Gamma_4$ consists of all subscripts $i$ with $p+1\le i \le p$ in ${\Pi}_{s_1s_2} \cup  {\Pi}_{t_1t_2}$ with  ${\rm val}(h_{i1}) > {\rm val}(h_{i2})$;

$\bullet$ $\Gamma_5$ consists of all subscripts $i$ with $p+1\le i \le p$ outside $\Gamma_3 \cup \Gamma_4$.

Define $K:= (\bigcup_{i \in \Gamma_1} J_i') \cup (\bigcup_{i \in \Gamma_2} J_i) \cup (\bigcup_{i \in \Gamma_3} L_i' \cup R_i) \cup 
(\bigcup_{i \in \Gamma_4} L_i \cup R_i')  \cup (\bigcup_{i \in \Gamma_5} L_i \cup R_i) \cup E_1$. It is routine to check that $K$ is an $(s_1,t_1;s_2,t_2)$-bicut in the 
network $(G,c)$. Observe that $c(K)$ equals $\sum_{i=1}^p {\rm val}(g_i)+ \sum_{i \in \Gamma_1} {\rm val}(h_i) + \sum_{i=p+1}^q [{\rm val}(g_{i1})+{\rm val}(g_{i2})]
+ \sum_{i \in \Gamma_3} {\rm val}(h_{i1}) + \sum_{i \in \Gamma_4} {\rm val}(h_{i2}) +c(E_1)$. Since ${\rm val}(h_i) = \min\{{\rm val}(h_{i1}), {\rm val}(h_{i2}) \}$
for $p+1\le i \le q$, we have $\sum_{i \in \Gamma_3} {\rm val}(h_{i1}) + \sum_{i \in \Gamma_4} {\rm val}(h_{i2})=\sum_{i \in \Pi_{s_1s_2}} {\rm val}(h_i)+
\sum_{i \in \Pi_{t_1t_2}} {\rm val}(h_i)$. It follows that $c(K)=\sum_{i=1}^p {\rm val}(g_i)+ \sum_{i=p+1}^q [{\rm val}(g_{i1})
+{\rm val}(g_{i2})] +w(E_1)=\sum_{i=1}^p {\rm val}(g_i)+ \sum_{i=p+1}^q [{\rm val}(g_{i1}) +{\rm val}(g_{i2})] + {\rm val}(l_1) + {\rm val}(l_2)
=\sum_{i=1}^p {\rm val}(g_i)+ \sum_{i=p+1}^q [{\rm val}(g_{i1}) +{\rm val}(g_{i2})] + {\rm val}(l_1^*) + {\rm val}(l_2^*)= {\rm val}(z_1) + {\rm val}(z_2)$.
Therefore, $(z_1,z_2)$ is a maximum integral $(s_1,t_1;s_2,t_2)$-biflow and $K$ is a minimum $(s_1,t_1;s_2,t_2)$-bicut in the network $(G,c)$. 

For each $1\le i \le p$, we need to apply Lemma \ref{triflow} once to $B_i$. For each $p+1 \le i \le q$, we need to apply Lemma \ref{triflow} twice to 
$B_i$. The problems of constructing $h_i$ (see the paragraph below (7)) and $l^*$ (see the paragraph above (9)) can both be reduced to maximum flow problems. 
For $1\le i \le q$, let $n_i$ and $m_i$ denote the numbers of vertices and edges in $B_i$, respectively. Then the whole algorithm runs in $O(\sum_{i=1}^q n_im_i)$ 
time. Note that $n=4+\sum_{i=1}^p (n_i-3)+ \sum_{i=p+1}^q (n_i-4)$.  So $\sum_{i=1}^q n_i \le n+4q \le 5n$. Hence the running time of our algorithm is $O(nm)$. \qed

\vskip 4mm
In Theorem \ref{Thm3} (and hence in Algorithms I, II and III), we assume that $G=(V,E)$ is a graph with two specified nonadjacent edges $e_1=s_1t_1$ 
and $e_2=s_2t_2$, and that $G$ is 2-connected and every 2-separation of $G$, if any, separates $e_1$ from $e_2$. If $e_i=s_it_i$ is not an edge of
the input graph addressed in Theorem \ref{Thm2} for $i=1$ or $2$, we may simply add it and define $c(e_i)=0$; the resulting graph remains to be a Seymour 
graph and is also $K_4^*$-free. To complete the proof, we still need to show that indeed the connectivity of the input graph can be lifted to meet the conditions. 

\vskip 2mm
{\bf Proof of Theorem \ref{Thm2}}. Let $G=(V,E)$ be a Seymour graph with two sources $s_1, s_2$ and two sinks $t_1,t_2$ and let $c:\, E \rightarrow \mathbb Z_+$ 
be a capacity function. To find a maximum integral $(s_1,t_1;s_2,t_2)$-biflow in the network $(G,c)$, clearly we may assume that $G$ is connected. 

Suppose $(G_1, G_2)$ is a $1$-separation of $G$. Let $c_i$ the restriction of $c$ to $G_i$ for $i=1,2$. If $G_1$ or $G_2$, say the former, contains neither 
$e_1$ nor $e_2$, then the maximum integral biflow problem on $(G,c)$ can be reduced to that on $(G_2, c_2)$. As the separation can be found in $O(n+m)$ time 
by Tarjan's algorithm \cite{T}, our algorithm requires $O(n+m)+O(n_2^3+m_2^2)=O(n^3+m^2)$ time by induction, where $n_2$ and $m_2$ are the numbers of vertices
and edges of $G_2$, respectively. So we may assume that $e_i \in E(B_i)$ for $i=1,2$. Thus to find a maximum integral biflow on $(G,c)$, it suffices to find 
a maximum integral $s_i$-$t_i$ flow $f_i$ in $(G_i, c_i)$ for $i=1,2$. Therefore, we may assume that $G$ is $2$-connected.    

Suppose $(G_1, G_2)$ is a $2$-separation of $G$ that does not separate $e_1$ for $e_2$, say neither is contained in $G_2$.  Let $\{s_3, t_3\}$ be the
corresponding $2$-separating set. Once again, let $c_i$ the restriction of $c$ to $G_i$ for $i=1,2$. We first find a maximum integral $s_3$-$t_3$ flow $g$ in
the network $(G_2,c_2)$. Let $G_1'$ be the union of $G_1$ and edge $s_3t_3$, and let $c_1'$ be the capacity function obtained from $c_1$ by defining 
$c_1'(s_3t_3):=c_1(s_3t_3)+{\rm val}(g)$, where $c_1(s_3t_3)=0$ if $s_3$ and $t_3$ are nonadjacent in $G$. We then find a maximum integral biflow $(f_1', f_2')$
in the network $(G_1', c_1')$, and finally get an integral biflow $(f_1, f_2)$ in the network $(G, c)$ by replacing certain units of arc flow $f_i'(s_3t_3)$ 
with a subflow of $g$ so that ${\rm val}(f_i)={\rm val}(f_i')$ for $i=1,2$. Obviously, $(f_1, f_2)$ is a maximum integral biflow in the network 
$(G, c)$. Thus the maximum integral biflow problem on $(G,c)$ can be reduced to that on $(G_1',c_1')$. Let $n_i$ and $m_i$ denote the numbers of vertices
and edges in $G_i$, respectively, for $i=1,2$. By induction, it takes $O(n_1^3+(m_1+1)^2)$ to find the biflow $(f_1', f_2')$. Hence the whole algorithm
requires $O(n_1^3+(m_1+1)^2)+O(n_2m_2)=O(n^3+m^2)$ time. Therefore, we may assume that every $2$-separation of $G$ separates $e_1$ for $e_2$.

By Theorem \ref{Thm3}, it takes $O(n^3)$ time to exhibit the global structure as described in (i)-(iv). We can then apply Algorithm I, II or III to 
the network $(G, c)$ to find a maximum integral biflow, which requires $O(m^2)$ time.  So the running time of the whole algorithm is $O(n^3+m^2)$. \qed

\end{document}